\documentclass[11pt,reqno]{amsart}
\usepackage{amsmath,amssymb,mathrsfs,graphicx}
\usepackage[utf8x]{inputenc}
\usepackage[T1]{fontenc}
\usepackage[USenglish,english]{babel}
\usepackage{amsmath}
\usepackage{amsfonts}
\usepackage{latexsym}
\usepackage{cite}
\usepackage{float}
\usepackage{tikz}
\usepackage{pgfplots}
\usepackage{caption}
\usepackage{booktabs} 
\usepackage{subcaption}
\usepackage{mathtools}      
\mathtoolsset{showonlyrefs} 
\usepackage[export]{adjustbox}
\usepackage{amsthm}
\usepackage{amssymb,amscd}
\usepackage{mathtools}
\usepackage{xargs}
\usepackage{graphicx}
\usepackage{color}
\usepackage{enumitem}
\usepackage{multirow}
\usepackage{lmodern}

\usepackage[text={6.5in, 8in}, centering]{geometry}
 
\usepackage{parskip}
\usepackage{microtype}

\usepackage{hyperref}

\newtheorem*{thm*}{Theorem}

\usepackage{appendix}

\theoremstyle{definition}

\newcommand{\A}{\mathbf A}
\newcommand{\Aa}{\mathbf a}
\newcommand{\bb}{\mathbf b}
\newcommand{\M}{\mathbf M}
\newcommand{\nn}{\mathbf n}
\newcommand{\T}{\mathbf T}
\newcommand{\Tss}{\mathbf{T}_{ss}}
\newcommand{\tcd}{T_f^{c,d}}
\newcommand{\Acd}{A_{c,d}}
\newcommand{\tpq}{t_{pq}}
\newcommand{\spq}{s_{pq}}
\newcommand{\pp}{\mathbf{e}_1}
\newcommand{\qq}{\mathbf{e}_2}
\newcommand{\Ts}{\T_s}
\newcommand{\Xt}{\mathbf X_t}
\newcommand{\X}{\mathbf X}
\newcommand{\Xs}{\mathbf X_s}
\newcommand{\Xss}{\mathbf X_{ss}}
\newcommand{\Tt}{\mathbf T_t}
\newcommand{\RR}{\mathbb R}
\newcommand{\Hbb}{\mathbb H}

\newcommand{\Talg}{\mathbf T_{alg}}

\newcommand{\Xalg}{\mathbf X_{alg}}

\DeclareMathOperator{\arccosh}{arccosh}
\DeclareMathOperator{\arcsinh}{arcsinh}
\DeclareMathOperator{\arctanh}{arctanh}

\title[The Schr\"{o}dinger map for helical polygons in the hyperbolic space]{On the Schr\"{o}dinger map for regular helical polygons in the hyperbolic space}
\author[S. Kumar]{Sandeep Kumar}
\address[S. Kumar]{BCAM - Basque Center for Applied Mathematics, Alameda de Mazarredo 14, 48009 Bilbao, Spain.}
\email{skumar@bcamath.org}
\date{\today}	
\begin{document}
\newenvironment{red}{\textcolor{red}}

\maketitle
\begin{abstract}

The main purpose of this article is to understand the evolution of $\Xt = \Xs \wedge_- \Xss,$ with $\X(s,0)$ a regular polygonal curve with a nonzero torsion in the 3-dimensional Minkowski space. Unlike in the case of the Euclidean space, a nonzero torsion now implies two different helical curves. This generalizes recent works by the author with de la Hoz and Vega on helical polygons in the Euclidean space as well as planar polygons in the Minkowski space. 

Numerical experiments in this article show that the trajectory of the point $\X(0,t)$ exhibits new variants of Riemann's non-differentiable function whose structure depends on the initial torsion in the problem. As a result, we observe that the smooth solutions (helices, straight line) in the Minkowski space show the same instability as displayed by their Euclidean counterparts and curves with zero-torsion. These numerical observations are in agreement with some recent theoretical results obtained by Banica and Vega.

\end{abstract}

\section{Introduction}

We consider the flow
\begin{equation}
\label{eq:VFE}
\Xt = \Xs \wedge_{+} \Xss,
\end{equation}
where $\X(s,t): \RR\times\RR \to \RR^3$ is an arc-length parameterized curve representing a vortex filament in a real fluid, i.e., smoke rings, helical vortices, and tornadoes, etc. Here, $s$ is the arc-length parameter, $t$ time, $\wedge_{+}$  the usual cross-product in the Euclidean space. Equation \eqref{eq:VFE} is called the vortex filament equation (VFE) or local induction approximation (LIA). Its first appearance dates back to 1906 where it was derived as a local approximation of the Biot--Savart integral describing the dynamics of a vortex filament under Euler equations\cite{darios}. 

With the tangent vector $\T = \Xs$, differentiating \eqref{eq:VFE} yields
\begin{equation}
\label{eq:SMP}
\Tt = \T \wedge_{+} \Tss,
\end{equation}
known as the Schr\"{o}dinger map equation onto the sphere. During the time evolution, the magnitude of $\T$ remains conserved and we can assume that it takes values on the unit sphere. Equation \eqref{eq:SMP} is a special case of Landau--Lifshitz equation in ferromagnetism \cite{landau}, and by writing it in a more geometric way:
\begin{equation}
\Tt = \mathbf{J} \mathbf{D}_s \Ts,
\end{equation}
where $\mathbf{J}$ is the complex structure of the sphere, $\mathbf{D}_s$ is the covariant derivative, it can be extended to other settings \cite{Khesin}. For instance, when the target space is chosen as $\Hbb^2=\{ (x_1,x_2,x_3):-x_1^2+x_2^2+x_3^2=-1, x_1>0  \}$, i.e., a unit sphere in the Minkowski 3-space $\RR^{1,2} = \{(x_1,x_2,x_3): ds^2 = -dx_1^2 + dx_2^2 + dx_3^2\}$, the equivalent of \eqref{eq:SMP} is \cite{Ding1998}
\begin{equation}
\label{eq:SMP-hyp}
\Tt = \T \wedge_{-} \Tss, 
\end{equation}
where $\T\in\Hbb^2$, and corresponding $\X\in\RR^{1,2}$ solves
\begin{equation}
\label{eq:VFE-hyp}
\Xt = \Xs \wedge_{-} \Xss.
\end{equation}
Here, the Minkowski cross product $\wedge_{-}$ is
$$
\Aa \wedge_- \bb = ( - (a_2b_3 - a_3b_2), a_3b_1 - a_1b_3, a_1b_2-a_2b_1), \ \Aa, \bb \in \RR^{1,2}. 
$$
The Minkowski pseudo-scalar product is given by $\Aa \circ_- \bb = -a_1b_1 + a_2 b_2 + a_3b_3,$
which defines 
$| \Aa |_0^2 = \Aa \circ_{-} \Aa$.
Let us recall that a vector in $\RR^{1,2}$ is \textit{space-like} (\textit{time-like}) if $|\cdot |_0^2$ is positive (negative), or, \textit{light-like} if $|\cdot|_0^2$ is zero. A time-like vector, in turn, can be further classified as \textit{positive (negative) time-like} depending on the sign of its first component. For example, the definition of $\Hbb^2$ involves the positive time-like vectors, and since in this article we work only with them, for the ease of notation, they are referred to as time-like. Moreover, as $\T\in\Hbb^2$, the corresponding curve $\X$ is termed as a time-like curve. Then, if $\X$ has a non-vanishing curvature $\kappa$, torsion $\tau$, the space-like normal vector $\nn$, binormal vector $\bb$ and $\T$ form an orthogonal system  and satisfy \cite{Lopez} 
\begin{equation}
\label{mat:SF-chap-hyp}
\begin{pmatrix} \T \\ \nn \\ \bb \end{pmatrix}_s = 
\begin{pmatrix} 0 & \kappa & 0\\ \kappa & 0 & \tau \\  0& -\tau & 0 \end{pmatrix} 
{.}\begin{pmatrix} \T \\ \nn \\ \bb \end{pmatrix}.
\end{equation}
Furthermore, with the Hasimoto transformation
\begin{equation}
\label{eq:fila-fun-hyp}
\psi(s,t) = \kappa(s,t) e^{i\int_{0}^{s}\tau(s^\prime,t)ds^\prime},
\end{equation}
\eqref{eq:SMP-hyp}--\eqref{eq:VFE-hyp} are related to the defocusing nonlinear Schr\"{o}dinger (NLS) equation \cite{hasimoto}:
\begin{align}\label{eq:NLS-hyp}
\psi_t = i\psi_{ss} - \frac{i}{2} \psi (|\psi|^2+A(t)), \ A(t)\in\RR.
\end{align} 
From \eqref{eq:fila-fun-hyp}, one can note that some explicit solutions of \eqref{eq:VFE}--\eqref{eq:VFE-hyp} correspond to the curves with $\kappa=a$, $\tau=b$, where $a$, $b \in\RR$. Unlike in the Euclidean case, when both $a$ and $b$ are nonzero, the corresponding curve $\X\in\RR^{1,2}$ can be categorized in two different ways. In other words, due to the symmetries of $\Hbb^2$, if the axis of rotation of the helical curve is a time-like vector then it is called a \textit{circular helix} and a \textit{hyperbolic helix} if it is a space-like vector \cite{Lopez}. 


Another important class is the one-parameter family of the self-similar solutions which are characterized by the parameter $c_0$. Recall that \eqref{eq:VFE}, \eqref{eq:VFE-hyp} are time-reversible, i.e., if $\X(s,t)$ is a solution then so is $\X(-s,-t)$. As a result, for both Euclidean and hyperbolic settings, it has been shown that as the time $t$ tends to zero, the solution curve develops a corner and turns into two non-parallel straight lines meeting at $s=0$ \cite{delahoz2007}; in this work, we refer to this as \textit{one-corner problem}. This implies that at $t=0$, the corresponding tangent vector is a Heaviside-type function and the initial solution of the NLS equation will have a Dirac delta located at $s=0$. The problem has been well-studied both from a theoretical and numerical point of view \cite{DelahozGarciaCerveraVega09,BV3}, and one of the most important outcomes is the relationship between the angle $\theta$ formed at the corner and the parameter $c_0$ which in the hyperbolic space is given by
\begin{equation}
\label{eq:1-corner-exp}
\cosh(\tfrac{\theta}{2}) = e^{-\pi c_0^2/2}.
\end{equation}
The motivation of this article stems from the recent work on curves with multiple corners which has revealed several interesting properties of \eqref{eq:VFE}--\eqref{eq:VFE-hyp} \cite{HozVega2014,HozKumarVega2019,HozKumarVega2020}. For instance, in \cite{HozKumarVega2020}, following the ideas of the Euclidean case \cite{HozVega2014}, the evolution of \eqref{eq:SMP-hyp}--\eqref{eq:VFE-hyp} was addressed for an initial datum being a regular planar polygon in $\mathbb{R}^{1,2}$, i.e., \textit{planar $l$-polygon} (from now on, referred to as \textit{planar polygon problem}). Characterized by the parameter $l>0$, unlike its Euclidean counterpart, a planar $l$-polygon is open \cite{FF}. At the level of the NLS equation, the corresponding initial data is a Dirac comb with constant coefficients, as a result, using the Galilean invariance and using some analytical-algebraic techniques the evolution at the rational times can be obtained up to a rigid movement. Numerical experiments confirm that as the polygonal curve evolves it develops more sides at each rational multiple of the time period, a behaviour reminiscent of the periodic Talbot effect in optics \cite{ET}. 
 
This complex dynamics is also present in the evolution of a point located on the polygonal curve in both Euclidean and hyperbolic cases \cite{HozVega2014,HozKumarVega2019,HozKumarVega2020}. For example, in the planar polygon problem, the curve $\X(0,t)$ was compared with the complex version of Riemann's non-differentiable function \cite{Du}
\begin{equation}
\phi(t) = \sum_{k=1}^{\infty} \frac{e^{\pi i k^2 t}}{i\pi k^2}, \ t\in[0,2].
\end{equation}
\noindent With adequate numerical simulations, a strong evidence was given that as the parameter $l$ approaches zero, the trajectory converges to $\phi(t)$. Note that $\phi(t)$ is an analytical object and was proved multifractal in \cite{Ja}, see also \cite{DE1} for some recent results on its regularity and geometric properties. Moreover, new variants of $\phi(t)$ have been discovered when the initial datum is a helical polygon in $\mathbb{R}^3$ \cite{HozKumarVega2019}. In fact, a strong numerical evidence is given that the structure of these new curves depends on the number of sides or the torsion of the initial curve. On the other hand, due to the lack of regularity mixed with periodicity, at a theoretical level, the problem is very challenging. In a recent work, for initial data as polygonal lines that at any given time are asymptotically close at infinity to two straight lines, by choosing an appropriate solution space, the well-posedness of the problem has been proved; moreover, it is  showed that the corresponding self-similar solutions have finite  renormalized energy \cite{BanicaVega2018,BanicaVega2019}. These non-trivial results are indeed useful as the principle of energy established there, allows us to compute the angle between any two sides of the polygonal curve for any rational time. Furthermore, very recently, building upon these ideas, by approximating the smooth curves, i.e., circle, straight line, and helix with non-closed polygonal lines, the existence of $\phi(t)$ and its invariants have been proved using  theoretical arguments rigorously \cite{BanicaVega2020}. Let us also mention that the presence of $\phi(t)$ was also observed in the linear momentum of the regular planar polygon in \cite{HozVega2018} (see \cite[Section 4.6]{KumarPhd} for the hyperbolic case). In the same work, authors also established a connection between the one-corner and multiple corner problems and utilized this to obtain some analytical results.

Thus, bearing in mind the theoretical results, in both Euclidean and hyperbolic settings, the appearance of Riemann's function (or its variants) while approximating the smooth curves with polygonal ones, amounts us to say that the evolution of the former is not stable. In other words, given a particle located arbitrarily close to a smooth curve with zero, or non-zero torsion when measured in the right topology, its evolution always approaches to the graph of $\phi(t)$ (or its variants). Let us also mention that the non-zero torsion case discussed in \cite{HozKumarVega2019} has also been compared with the one of the Schr\"{o}dinger equation with quasi-periodic boundary conditions studied recently in \cite{BFP}.

In this article, we extend these results by considering the non-planar polygonal curves in the Minkowski space $\RR^{1,2}$. Using the previous notations, we term the two types of helical polygonal curves as a \textit{circular helical polygon (CHP)} and a \textit{hyperbolic helical polygon (HHP)}. For an arc-length parameterized CHP with $M$ sides, the tangent vector $\T$ are $2\pi$-periodic and take $M$ values on $\Hbb^2$. More precisely, the first component is a constant $b\in(1,\infty)$, the second and third components lie on a circle with a radius $a$ so that $a^2-b^2=-1$ (see Figure \ref{fig:H2}). On the other hand, when the parameter $b \in(0,\infty)$ corresponds to the second, or the third component of $\T\in\Hbb^2$, the remaining two components lie on a hyperbola with a radius $a$ so that $-a^2+b^2=-1$ (see Figure \ref{fig:H2}). The resulting curve $\X$ is a HHP whose side-length is denoted by a parameter $l>0$. Note that, in both cases, the parameter $b$ introduces torsion to the initial curve, for example, for a HHP, $b=0$ reduces back to the planar $l$-polygon case, for a CHP, $b=1$ to a straight line and for an intermediate value of $b$ the corresponding polygonal curve has a (hyperbolic/circular) helical shape. We describe the evolution of \eqref{eq:SMP-hyp}--\eqref{eq:VFE-hyp} for these curves as initial data and call it \textit{helical polygon problem}.

The structure of the article is the following. In Section \ref{sec:Sol-hel-pol}, we discuss the theoretical formulation of the problem, symmetries for both types of helical polygons and by working with the NLS equation, we construct the time evolution of the $\X$ and $\T$ at the rational times using algebraic techniques. In Section \ref{sec:Num-sol}, we compute the numerical evolution of the helical polygons and show a good agreement between the numerical and exact values of quantities such as the speed of the center of mass, the angle between any two sides of the polygon, etc. One of our main aims in this article, is to show that the smooth helical curves in $\RR^{1,2}$, exhibit the same instability as displayed by their Euclidean counterparts and curves with zero-torsion. We claim this for the helical curves in Section \ref{sec:Traj-X0t}, by showing the presence of new variants of Riemann's function in the evolution of a point located on the helical polygon. Section \ref{sec:Tstpq-qgg1} discusses the evolution of the tangent vector at rational time with a large denominator which can be seen as the approximation of irrational times. Finally, in Section \ref{sec:Num-rel}, we establish a numerical relationship between the helical polygon and one-corner problems and investigate the relevant consequences which also allows us to compare it with its equivalent in the Euclidean space.
 
\section{A solution of (\ref{eq:SMP-hyp})--(\ref{eq:VFE-hyp}) for Regular helical polygons in the Minkowski space}
\label{sec:Sol-hel-pol}

As the first objective of this article, we describe the dynamics of \eqref{eq:SMP-hyp}--\eqref{eq:VFE-hyp} when the initial datum is considered as a regular polygonal curve with a nonzero torsion. The evolution of these curves is a result of the Galilean symmetry in the set of solutions of \eqref{eq:NLS-hyp}, thus, this work can be seen as a natural extension of the zero-torsion case \cite{HozKumarVega2020}. Consequently, \cite[Thm. 1]{HozKumarVega2020} is also valid here with a change of \cite[(23)]{HozKumarVega2020} and \cite[(24)]{HozKumarVega2020} with \eqref{eq:psi-theta-stpq-hyp-hel} and \eqref{eq:rhoq-tpq-hel-hyp}, respectively. Next, we study the evolution of the two types of polygonal curves, after representing them in parametric form.

\subsection{Problem definition}
\label{sec:Pb-def-form}
\subsubsection{Hyperbolic helical polygon (HHP)}
An arc-length parameterized regular hyperbolic helical polygon is characterized by parameters $l$ and $b$ corresponding to its side-length and torsion, respectively. Due to the rotation invariance of $\wedge_-$ and symmetries of $\mathbb{H}^2$, without any loss of generality, we can choose $b$ as the third component of the tangent vector $\T(s,0)$ (see Figure \ref{fig:H2}). Thus, we write 
\begin{equation} \label{eq:T0-hyp-hel}
\T(s,0) = \left(a\cosh\left(\frac{l}{2}+s_k\right),a\sinh\left(\frac{l}{2}+s_k\right),b \right), \quad  \quad s_k<s<s_{k+1},
\end{equation}
where $s_k = k l$ for $k\in\mathbb{Z}$, $b>1,a^2-b^2=1$. 
The corresponding curve $\X(s,0)$ is a helical polygon with vertices located at
\begin{equation}\label{eq:X0-hyp-hel}
\X(s_k,0) =  \frac{l/2}{\sinh(l/2)} \left( a\sinh(s_k), a \cosh(s_k),b s_k \right),
\end{equation}
so that for any $s \in (s_k , s_{k+1})$, the corresponding point $\X(s, 0)$ lies on the line segment joining $\X(s_k,0)$ and $\X(s_{k+1},0)$. Here, we have taken \eqref{eq:T0-hyp-hel} and thus, \eqref{eq:X0-hyp-hel} in such a form so that the vertex located at $\X(0,0)$ lies on the $y$-axis. Note that $b\in(-\infty,\infty)$; however, in this work, we consider the case when $b>0$, since the case $b<0$ can be easily recovered using the symmetries. 

By denoting $\T_k = \T(s,0)$ for $s\in(s_k,s_{k+1})$, the time-like curvature angle $\rho_0$ between any two sides can be computed as 
\begin{align}\label{eq:rho0-hyp-hel}
	\cosh(\rho_0) &= -\T_k \circ_- \T_{k+1} 
	= 1+2a^2\sinh^2\left(\frac{l}{2}\right) \Longleftrightarrow  
	\rho_0 = 2 \arcsinh\left(a\sinh\left(\frac{l}{2}\right)\right).
\end{align}
The vector $\T_k \wedge_-\T_{k+1}$, is space-like for all $k$, and satisfy
\begin{equation}
\label{eq:Tk-Tk1-prop}
(\T_{k-1}\wedge_-\T_k) \circ_- (\T_k\wedge_-\T_{k+1}) < |\T_{k-1}\wedge_-\T_k |_0 \; |\T_k\wedge_-\T_{k+1} |_0, 
\end{equation} 
as a result, the angle between them is space-like \cite[Thm. 3.2.6 \& (3.2.6)]{HypGeoText1}, and we have the space-like torsion angle
\begin{equation} 
\label{eq:theta_def_hyp}
\theta_0 = \arccos\left(\frac{(\T_{k-1} \wedge_{-} \T_{k})\circ_- (\T_k \wedge_{-} \T_{k+1} )}{|\T_{k-1} \wedge_{-} \T_{k} |_0 \ |\T_k \wedge_{-} \T_{k+1} |_0}\right), 
\end{equation}
or,
\begin{equation} \label{eq:theta0-hyp-hel}
\theta_0 = 2 \arctan\left(b\tanh\left(\frac{l}{2}\right)\right). 
\end{equation}		
Moreover, $\rho_0$, $\theta_0$ do not depend on $k$ and satisfy
 $$
\cos\left(\frac{\theta_0}{2}\right)\cosh\left(\frac{\rho_0}{2}\right) = \cosh\left(\frac{l}{2}\right). 
$$	 
Like in the zero-torsion case, for a hyperbolic helical polygon, $\X$ and $\T$ are invariant under a hyperbolic rotation of angle $k l$ about the space-like $z$-axis (see \cite[Sec. 2.2.1]{HozKumarVega2020}). However, the third component $X_3$ has a translation symmetry, i.e., $X_3(s+kl,t) =  X_3(s,t) + k l b$, for all $k\in \mathbb{Z}$, $t\geq0$. 
\subsubsection{Circular helical polygon (CHP)}
\begin{figure}[htbp!] 
	\centering
	\includegraphics[width=0.51\textwidth, valign=t]{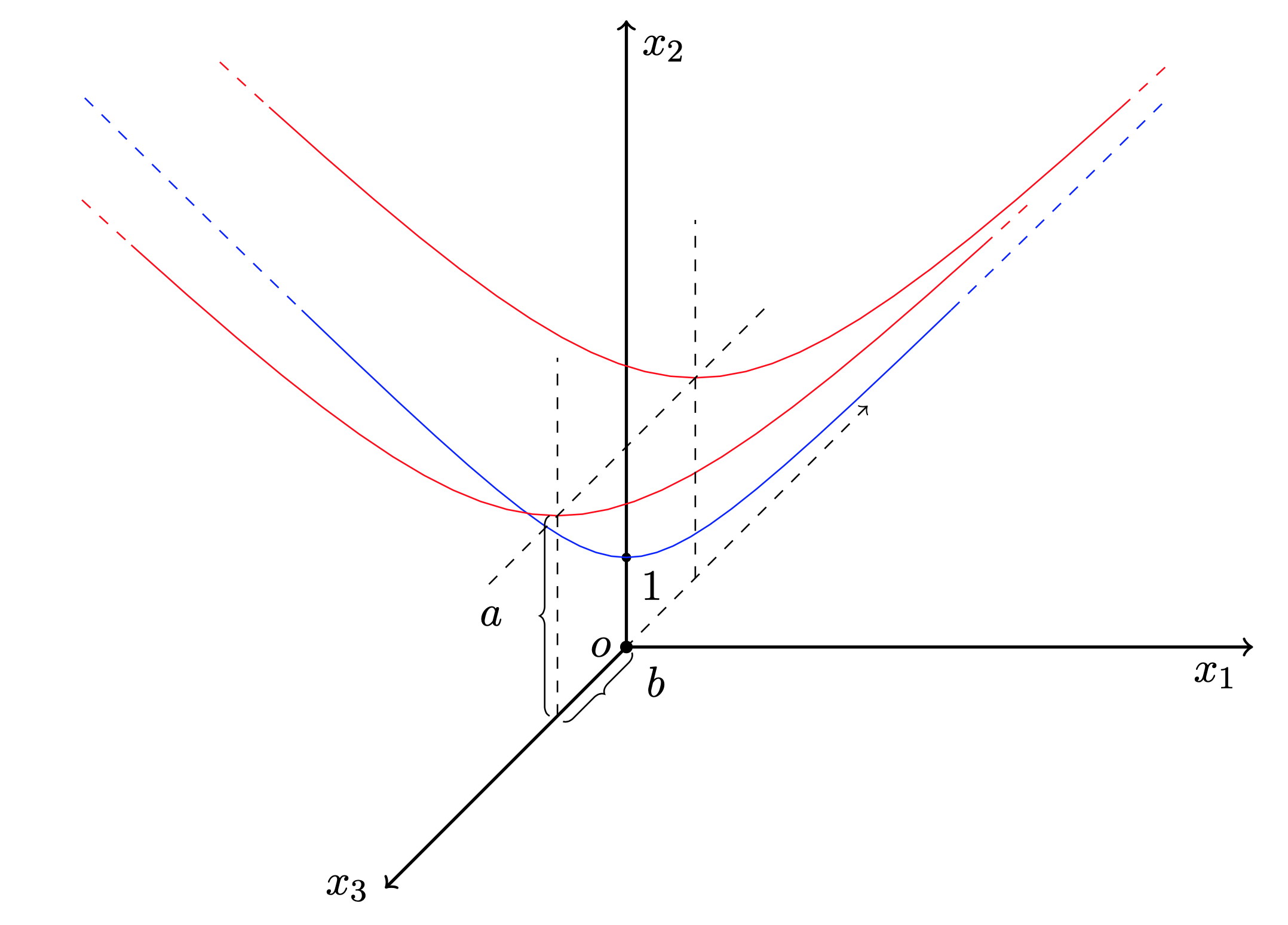}	
	\includegraphics[width=0.45\textwidth, valign=t]{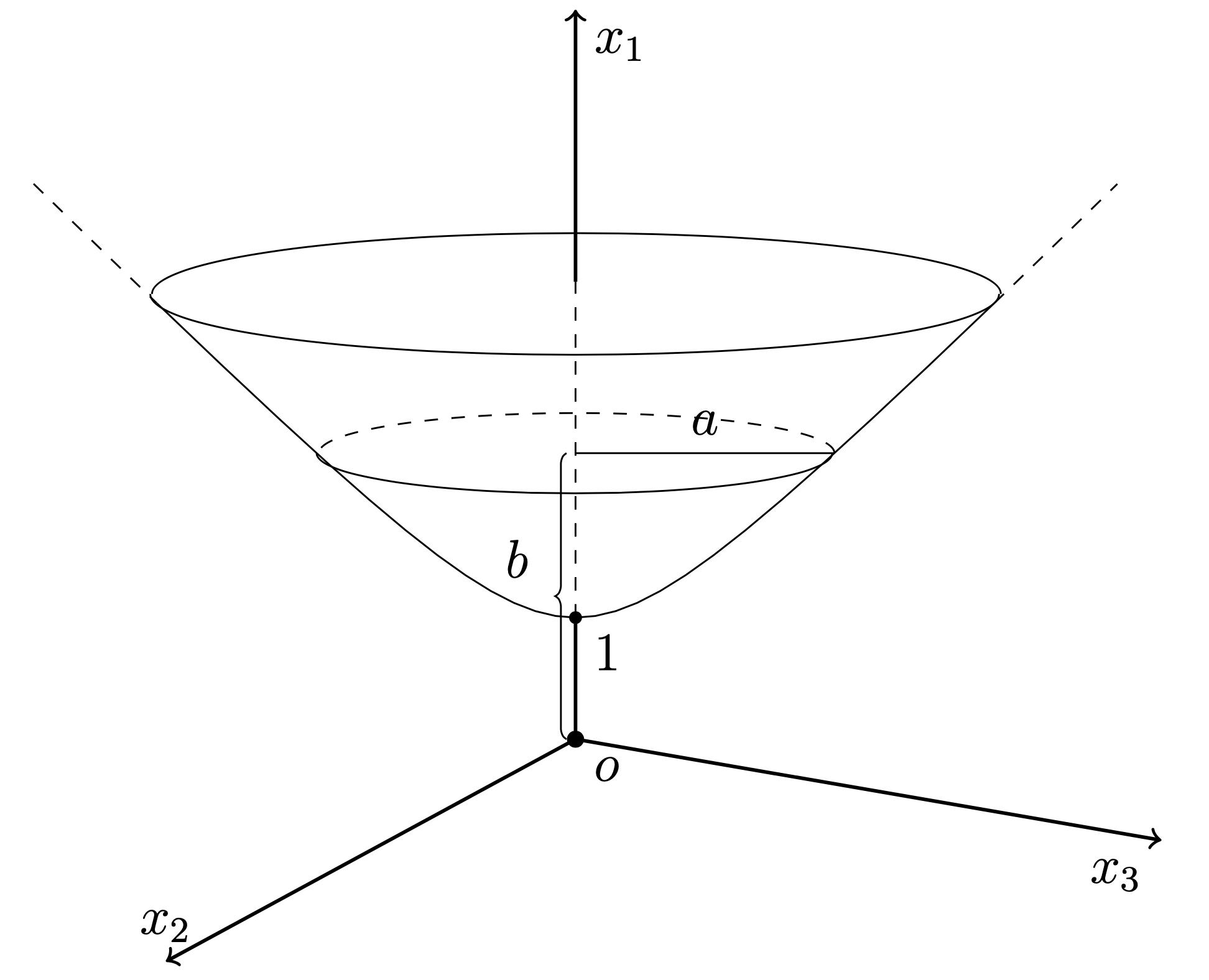}
	\caption{Left: Hyperbolic 2-space $\Hbb^2$ and the parameter $b\geq0$ such that $a^2-b^2=1$. The red (blue) unit hyperbolas are formed by the intersection of plane $z=\pm b$ ($z=0$) with $\Hbb^2$. Right:  $\Hbb^2$ and the parameter $b\geq1$ such that $b^2-a^2=1$. }
	\label{fig:H2}   
\end{figure}

We consider an arc-length parameterized regular circular helical polygon with $M$ sides and a torsion determined by the parameter $b$. In this case, $b$ corresponds to the first component of the $2\pi$-periodic tangent vector $\T(s,0)$ as in Figure \ref{fig:H2}, and we write 
\begin{equation} \label{eq:T0-circ-hel}
\T(s,0) = \left(b,a\cos\left(\frac{2\pi k}{M}\right),a\sin\left(\frac{2\pi k}{M}\right) \right), \quad  \quad s_k<s<s_{k+1},
\end{equation}
where $s_k = 2\pi k /M$ for $k=0,1,...,M-1$, $b>1,a^2-b^2=-1$. 
The corresponding curve $\X(s,0)$ is a helical polygon whose vertices are given by 
\begin{equation}\label{eq:X0-circ-hel}
\X(s_k,0) = \left(bs_k, \frac{a\pi\sin(\pi(2k-1)/M)}{M\sin(\pi/M)}, -\frac{a\pi\cos(\pi(2k-1)/M)}{M\sin(\pi/M)} \right),
\end{equation}
so that for any $s \in (s_k , s_{k+1})$, $\X(s, 0)$ lies on the line segment joining $\X(s_k,0)$ and $\X(s_{k+1},0)$. 

Since $\T\in\Hbb^2$, $|b|\geq1$; however, thanks to the symmetries, it is sufficient to work with only positive values of $b$. Thus, with $\T_k = \T(s,0)$, $s\in(s_k,s_{k+1})$, the time-like curvature angle $\rho_0$ between any two sides is
\begin{equation}
\label{eq:rho0-circ-hel}
\cosh(\rho_0) = -\T_k \circ_- \T_{k+1} 
= 1+2a^2\sin^2\left(\frac{\pi}{M}\right) \Longleftrightarrow   
\rho_0 = 2 \arcsinh\left(a\sin\left(\frac{\pi}{M}\right)\right). 
\end{equation}
As before, $\T_{k}\wedge_-\T_{k+1}$ is space-like for all $k$ and they satisfy \eqref{eq:Tk-Tk1-prop}, consequently, from \eqref{eq:theta_def_hyp}, the space-like torsion angle 
\begin{equation} \label{eq:theta_expression-circ-hel}
\theta_0 = 2 \arctan\left(b\tan\left(\frac{\pi}{M}\right)\right);
\end{equation}
and independent of $k$, both $\rho_0$, $\theta_0$ satisfy
$$
\cos\left(\frac{\theta_0}{2}\right)\cosh\left(\frac{\rho_0}{2}\right) = \cos\left(\frac{\pi}{M}\right). 
$$
Note that $\X$ and $\T$ are invariant under a rotation of an angle $2\pi k/M$ about the time-like $x$-axis and the first component $X_1$ has a translation symmetry, i.e., $X_1(s+\frac{2\pi k}{M},t) =   X_1(s,t) + \frac{2\pi k b}{M}$, for all $k\in\mathbb{Z}$, $t\geq0$. Due to the mirror invariance $\T(s,t)$ and $\T(-s,t)$ are symmetric about the XY-plane; consequently, $\X(s,t)$ and $\X(-s,t)$ are symmetric about the $z$-axis, for all $t$. An important corollary of these symmetries is that $\X(s+2\pi k,t) - \X(s,t) = 2 \pi b\, k (1,0,0)^T$, for all $k\in \mathbb{Z}$ which will be helpful later. 
	
Let us mention that we denote \eqref{eq:fila-fun-hyp} as $\psi_\theta$ when $\theta_0\neq0$ and $\psi$ for $\theta_0=0$. Furthermore, while working with polygonal curves, a more general form of Frenet--Serret formulas is preferred, i.e., \cite[(21)]{HozKumarVega2020} where the derivatives of the new unit space-like normal and binormal $\pp$, $\qq$, respectively, depend only on the unit time-like tangent vector $\T$ \cite{HozKumarVega2020}.  

\subsection{Problem formulation and the evolution at the rational times}
\label{sec:Psi-tpq-theta}
\subsubsection{Hyperbolic helical polygon}
Recall that in the zero-torsion case, the initial data for \eqref{eq:NLS-hyp} is given by 
\begin{equation}
\label{eq:psi0}
\psi(s,0) = c_0 \sum_{k=-\infty}^{\infty}\delta(s-lk), s\in \RR,
\end{equation}
and when $\theta_0\neq0$,
\begin{equation}
\label{eq:psi0-theta}
\psi_\theta(s,0) = c_{\theta,0}\; e^{i\gamma s}\sum_{k=-\infty}^{\infty}\delta(s-lk),
\end{equation}
where $\gamma = \theta_0 / l$, and $\theta_0$ is as in \eqref{eq:theta0-hyp-hel} and constants $c_0$, $c_{\theta,0}$ depend on the initial configuration of the corresponding polygonal curve. In particular, 
\begin{equation}
\label{eq:c0-ctheta0}
c_0 = \sqrt{\frac{2}{\pi}\ln\left(\cosh\left(\frac{l}{2}\right)\right)}, \ 
c_{\theta,0} = \sqrt{\frac{2}{\pi}\ln\left(\cosh\left(\frac{\rho_0}{2}\right)\right)},\end{equation}
where $\rho_0$ is as in \eqref{eq:rho0-hyp-hel}. Using \eqref{eq:psi0}--\eqref{eq:psi0-theta},
$\psi_\theta(s,0) = \frac{c_{\theta,0}}{c_0} e^{i\gamma s} \psi(s,0)$,
then, from the Galilean invariance of \eqref{eq:NLS-hyp}, we obtain 
\begin{equation} \label{eq:psi-theta-hyp-hel}
\psi_{\theta}(s,t) =\frac{c_{\theta,0}}{c_0} \ e^{i\gamma s- i \gamma^2 t} \psi(s-2\gamma\, t,t).
\end{equation}
Using \cite[(29)]{HozKumarVega2020}, 
\begin{equation} \label{eq:psi-theta-st-hyp-hel}
\psi_\theta(s,t) =   \frac{c_{\theta,0}}{c_0}  \hat{\psi}(0,t) e^{i\gamma s - i \gamma^2 t} \sum\limits_{k=-\infty}^{\infty}  e^{-i(rk)^2 t+irk(s-2\gamma t)},
\end{equation}	
where $r=2\pi/l$ and $\hat \psi(0,t)$ depends on time and because of the Gauge invariance of \eqref{eq:NLS-hyp}, can be taken real for all $t$ \cite{HozVega2014}.  	
Note that in \eqref{eq:psi-theta-hyp-hel}, $\psi$ on the right-hand side is $l$-periodic in space and $l^2/2\pi$-periodic in time, but due to the presence of $\gamma$ (or $\theta_0$), $\psi_{\theta}$ ($\X$, $\T$) are neither space nor time-periodic. However, as we will see in the numerical simulations, their structure repeats every time when the time period is increased by $l^2/2\pi$. Hence, with some abuse of notation, we call this important quantity the time-period and denote it by $T_f = l^2/2\pi$. 
	
Next, we compute $\psi_\theta(s,t)$ at rational multiples of $T_f$, i.e., at $t=\tpq = T_f (p/q)$, with $p\in\mathbb{Z}$, $q\in \mathbb{N}$, and $\gcd(p,q)=1$. Thus, substituting values of $\gamma$, $\tpq$, $\hat \psi(0,\tpq)=c_0/l$, and further simplification gives (see \cite{HozVega2014} for the intermediate steps)
\begin{align}\label{eq:psi-theta-stpq-hyp-hel}
\psi_\theta(s,t_{pq}) =
\begin{cases}
\frac{c_{\theta,0}}{\sqrt{q}} e^{i(\theta_0^2/(2\pi))(p/q)}  \sum\limits_{k=-\infty}^{\infty}
\sum\limits_{m=0}^{q-1}
\\
\quad e^{ i   (\xi_m+k\theta_0+m\theta_0/q)} \delta\left(s- \frac{l\theta_0 p}{\pi q} - l k -\frac{l m}{q}\right) &\text{if $q$ odd},
\\
\frac{c_{\theta,0}}{\sqrt{q/2}} e^{i(\theta_0^2/(2\pi))(p/(q/2))}  \sum\limits_{k=-\infty}^{\infty} \sum\limits_{m=0}^{q/2-1}
\\
\quad e^{ i   (\xi_{2m}+k\theta_0+2m\theta_0/q) }\delta\left(s-\frac{l\theta_0 p}{\pi q} - lk -\frac{2 l m}{q}\right)  &\text{if $\frac q2$ even},
\\
\frac{c_{\theta,0}}{\sqrt{q/2}} e^{i(\theta_0^2/(2\pi))(p/(q/2))} \sum\limits_{k=-\infty}^{\infty}\sum\limits_{m=0}^{q/2-1} 
\\
\quad e^{ i   (\xi_{2m+1}+k\theta_0+(2m+1)\theta_0/q) }\delta\left( s-\frac{l\theta_0 p}{\pi q} - lk -\frac{l  (2m+1)}{q}\right)  &\text{if $\frac q2$ odd}.
\end{cases}
\end{align}
It is important to mention here that the value of $\hat \psi(0,\tpq)$ has been obtained in \cite{HozKumarVega2020} which in turn follows from a conservation law established for polygonal lines in \cite{BanicaVega2018,BanicaVega2019}. Moreover, due to the Galilean transformation, at time $\tpq$, a corner initially located at $lk$, $k\in\mathbb{Z}$ is translated by $\spq=l\theta_0p/\pi q$. We refer to this extra movement as \textit{Galilean shift}. 
	 
On the other hand, from the expression \eqref{eq:psi-theta-stpq-hyp-hel}, one can note that at any rational time $\tpq$, the modulus of the coefficients of newly formed Dirac deltas is a constant, whic{{h we denote by $c_{\theta,q}$. More precisely, 
\begin{equation}
\label{eq:c0-tpq-hel-hyp}
c_{\theta,q}=
\begin{cases}
\frac{c_{\theta,0}}{\sqrt{q}} \ & \text{if $q$ odd,} \\
\frac{c_{\theta,0}}{\sqrt{q/2}} \ & \text{if $q$ even.}
\end{cases}
\end{equation}
From the discussion in \cite[Section 2.3.1]{HozKumarVega2019}, the choice of \eqref{eq:c0-ctheta0} follows from \eqref{eq:1-corner-exp} which also holds true whenever there is a corner formation, for instance, in our case, at any rational time $\tpq$. Then, denoting the corresponding angle between any two sides at those times by $\rho_q$, we have 
\begin{equation}
\label{eq:c0-rhoq-tpq-hel-hyp}
c_{\theta,q} = \sqrt{\frac{2}{\pi} \log\left(\cosh\left(\frac{\rho_q}{2}\right)\right)},
\end{equation}
which upon combining with \eqref{eq:c0-tpq-hel-hyp} and \eqref{eq:c0-ctheta0} yields
\begin{equation}
\label{eq:rhoq-tpq-hel-hyp}
\rho_q=
\begin{cases}
\arccosh\left(\cosh^{1/q}\left(\rho_0/2\right)\right) \ & \text{if $q$ odd,} \\
\arccosh\left(\cosh^{2/q}\left(\rho_0/2\right)\right) \ & \text{if $q$ even.}
\end{cases}
\end{equation}
\subsubsection{Circular helical polygon}
\label{sec:Pb-f-CHP}
Continuing with the same notations, we write the initial data for a CHP with $M$ sides, as 
$$
\psi_\theta(s,0) = c_{\theta,0} e^{i\gamma s}\sum_{k=-\infty}^{\infty}\delta(s-2 \pi k/M), \ s\in[0,2\pi),
$$
where $\gamma = M \theta_0 / 2\pi$, and $\rho_0$, $\theta_0$, $c_{\theta,0}$ are as in \eqref{eq:rho0-circ-hel}, \eqref{eq:theta_expression-circ-hel}, \eqref{eq:c0-ctheta0}, respectively. Then, following the same as steps as above, we note that when $\gamma\in(0,1)$, the expression for $\psi_\theta(s,t)$ is neither periodic in space nor time; however, when $\gamma=1$, it is periodic both in space and time. Furthermore, the quantity $T_f=2\pi/M^2$, and evaluating it at the rational times $t=\tpq$, gives an expression for $\psi_\theta(s,\tpq)$ as in \eqref{eq:psi-theta-stpq-hyp-hel} which implies that at those times, the equally-spaced $M$ Dirac deltas in $s\in[0,2\pi)$ at $t=0$ turn into equally-spaced $Mq$ Dirac deltas (if $q$ odd) or $Mq/2$ Dirac deltas (if $q$ even). Moreover, the Galilean shift $\spq=2\theta_0p/Mq$, and the absolute value of corresponding coefficients of the Dirac deltas is constant which determines the angle $\rho_q$ as in \eqref{eq:rhoq-tpq-hel-hyp} with $\rho_0$ as in \eqref{eq:rho0-circ-hel}. 
\subsection{Algebraic solution}
\label{sec:Alg-sol}
Following the approach in \cite{HozKumarVega2019,HozKumarVega2020}, we construct the evolution using algebraic techniques and denote it by $\Xalg$, $\Talg$. In this direction, given any rational time $\tpq$ with $q$ odd (similarly, for $q$ even), we first define 
\begin{equation}\label{eq:Scaling-hel-hyp}
\Psi_\theta  (s,\tpq) = \frac{\rho_q}{c_{\theta,q}} e^{-i(\theta_0^2/(2\pi))(p/q)}\psi_\theta(s,\tpq).
\end{equation}
Since $ \lim_{q\to \infty} \rho_qe^{-i(\theta_0^2/(2\pi))(p/q)}/c_{\theta,q}<\infty$, $\Psi_\theta$ is well-defined and \eqref{eq:Scaling-hel-hyp} is equally  valid for both CHP and HHP cases. 

Using \eqref{eq:Scaling-hel-hyp} and expression for $\psi_\theta(s,\tpq)$, we can write $\Psi_\theta  (s,\tpq) =  (\alpha_{k,m} + i \beta_{k,m})(s,\tpq) = \rho_q e^{i\zeta_{k,m}}$. As a result, integrating the generalized Frenet--Serret frame gives the rotation matrix $\mathbf{R}_{k,m}$ (see \cite[(35)]{HozKumarVega2020}) which performs a rotation of time-like angle $\rho_q$ about a space-like rotation axis $(0,-\sin (\zeta_{k,m})$, $\cos (\zeta_{k,m}))$. In other words, it describes the transition across a corner at $s=  \tilde s_k$ where for the CHP, $\tilde s_k = (2\pi  (k+1) /Mq + \spq)^-$; on the other hand, for the numerical computations we work with a truncated HHP with $M$ sides, length $L=lM$, and $\tilde s_k = l (k+1) /q + \spq$, for $k=0,1,\ldots,Mq-1$. Thus, by taking the basis vectors $\tilde{\T}(s)$, $\tilde{\mathbf{e}}_1(s)$, $\tilde{\mathbf{e}}_2(s)$, as the identity matrix at $s = \spq^-$, and we obtain their values at $s=  \tilde s_k^-$, $k=0,1,\ldots,Mq-1$, by the action of $Mq$ rotation matrices. 

Furthermore, the vertices of $\tilde{\X}$ (i.e., $\X$ up to a rigid movement) are shifted by the Galilean shift $\spq$, and can be obtained by
\begin{align}
\label{eq:Xalg-from-Talg}
\begin{cases}
\tilde{\X}(\spq) = \tilde{\X}(0) + \spq \tilde{\T}(\spq^-), \\
\tilde{\X}(\tilde s_k) = 	\tilde{\X}(\tilde s_{k-1}) + \Delta s \tilde{\T}(\tilde s_k^-), & k=0,1,\ldots,Mq-1,
\end{cases}
\end{align}
with $\tilde{\X}(0)=(0,0,0)$, $\Delta s = 2\pi/Mq$ for a CHP and $\Delta s = l / q$ for a HHP. In this way, we get the positions of the vertices corresponding to the interval $s\in[\spq,2\pi+\spq]$ for the CHP and $s\in[-L/2+\spq,L/2+\spq]$ for the HHP. Next, by calculating the non-vertex points through the linear interpolation, and using the relevant symmetries mentioned before, the corresponding $\tilde{\X}(s)$ can be obtained for all $s\in[0,2\pi]$ and $s\in[-L/2,L/2]$.
			
Let us mention that the computation of $\tilde{\X}(2\pi+\spq)$ allows us to define the correct rotation matrix $\mathbf{L}_1$ so that the CHP can be aligned in such a way that its axis of rotation is parallel to the $x$-axis. To be precise, $\mathbf{L}_1$ is a rotation of an angle equal to the one between $\mathbf{v}= \tilde{\X}(2\pi+\spq)-\tilde{\X}(\spq)$ and $(1,0,0)^T$, about an axis orthogonal to a plane spanned by these two vectors. Consequently, $\T = \mathbf{L}_1 \cdot \tilde{\T}$, $\X = \mathbf{L}_1 \cdot \tilde{\X}$; then, after adding a constant to the second and third components of $\X$, so that the means of these components over a spatial period are zero, we obtain $\X$ and $\T$, up to a vertical movement and a rotation about the $x$-axis, where the former can be computed at any time $\tpq$ from the speed of the center of mass $c_M$. Similarly, for the HHP case, we align the curve $\tilde \X$ so that its first two components lie in the plane orthogonal to $z$-axis. In this direction, we compute the corresponding rotation $\mathbf L_1$ by bearing in mind that at any time $\tpq$ the space-like vectors $\tilde \T(-L/2)-\tilde \T(-l)$ and $\tilde \T(L/2)-\tilde \T(l)$ lie in the plane orthogonal to the $z$-axis. Thus, after adding the movement of the center of mass, we obtain the curve $\X = \mathbf L_1 \cdot \tilde \X$ and $\T = \mathbf L_1 \cdot \tilde \T$, up to a rotation about the $z$-axis. 

			
\section{Numerical experiments}
\label{sec:Num-sol}

Due to the geometrical nature of the two helical polygons that are taken as the initial data, \eqref{eq:SMP-hyp}--\eqref{eq:VFE-hyp} are approximated numerically through different methods. There have been several articles dedicated to the numerical treatment of \eqref{eq:VFE}--\eqref{eq:VFE-hyp}, for instance, see \cite{DelahozGarciaCerveraVega09,HozKumarVega2017}; however, the methods proposed in \cite{HozVega2014,HozKumarVega2020} are the most suitable ones for our purpose. Thus, as the tangent vector $\T$ for a CHP is $2\pi$-periodic, we can employ a pseudo-spectral space discretization and a fourth-order Runge--Kutta method in time for initial data \eqref{eq:T0-circ-hel}--\eqref{eq:X0-circ-hel}. In other words, the space interval $[0,2\pi]$ is divided into $N$ equally-spaced nodes, $s_k = 2\pi k/N$, $k=0, 1, \ldots, N-1$, and the time interval $[0,T_f]$, with $T_f = 2\pi/M^2$ into $N_t+1$ different time steps $t_n = n T_f/N_t$, $n=0,1, \ldots, N_t$. Thanks to the symmetries of $\T$, the approximation of the space derivatives of $\T$ through the discrete Fourier transform can be reduced from $N$ to $N/M$ elements \cite{HozVega2014} which in turn can be implemented with the \texttt{fft} algorithm in MATLAB. As a result, the computational cost can be reduced quite effectively. 

On the other hand, in the case of a HHP, given the parameters $l$ and $b$, a truncated curve $\X$ with $M$ sides so that its length $L=lM$, is considered. As proposed in \cite{HozKumarVega2020}, we use a finite difference space discretization and a fourth-order Runge--Kutta method in time with fixed boundary conditions on $\T$. Thus, the space interval $[-L/2,L/2]$ is divided into $N+1$ equally-spaced nodes $s_k = -L/2+kL/N$, $k=0,1,\ldots,N$, and with $T_f = l^2/2\pi$, as before $[0,T_f]$ is divided into $N_t$ equal length intervals. The $N+1$ values of the initial data \eqref{eq:T0-hyp-hel}--\eqref{eq:X0-hyp-hel} have been calculated as in \cite[Section 3.1]{HozKumarVega2020}. For more accurate results, a large value of $M$ is desired and for the symmetries mentioned in Section \ref{sec:Pb-def-form} we have taken it to be an even number.   

In both cases, at the end of each time step $t_n$, $\T$ is renormalized to make sure that it lies in $\Hbb^2$. Furthermore, the stability constraints also do not change from those in \cite{HozVega2014} and \cite{HozKumarVega2020}, respectively. Finally, unlike in the Euclidean case, due to the exponential growth of the Euclidean norm of $\T$, working with all values of the parameter $b$ (and $l$ for a HHP) is not possible numerically. Recall that similar restrictions were also observed in \cite{DelahozGarciaCerveraVega09,HozKumarVega2020}.  

Note that a CHP is characterized by parameters $M$ and $b$, so for a fixed $b$, as $M$ increases, the resulting curve tends to a circular helix. On the other hand, for a fixed $M$, when $b$ approaches $1$, it tends to a straight line. Similarly, in the case of a HHP, for a given $b$, as parameter $l$ goes to zero, we obtain a hyperbolic helix and with $b=0$, the planar $l$-polygon is recovered \cite{HozKumarVega2020}. In our numerical simulations, we have analyzed each of the three cases by calculating the relevant quantities. Thus, besides the fact that at any rational times $\tpq$, depending on the parity of $q$, the $M$-sided initial curve develops $Mq$ or $Mq/2$ number of sides, we observe that the evolution in both type of curves also experiences lack in space and time periodicity as commented in the previous section. This in turn leads to the terms \textit{Galilean shift} and \textit{phase shift}. While the former corresponds to the translation of a corner initially located at $s=0$ along the curve, the latter stands for its angular movement in the XY-plane for a HHP and in the YZ-plane for a CHP. This has been depicted in Figure \ref{fig:Gal-phase-HHP} (the left-hand side) for parameters $M=48$, $l=0.4$, $\theta_0=\pi/4$ for a HHP (only inner vertices) and for a CHP (the right-hand side) with parameters $M=3$, $b=1.2$.  
\begin{figure}[htbp!] 
\centering
\includegraphics[width=0.43\textwidth]{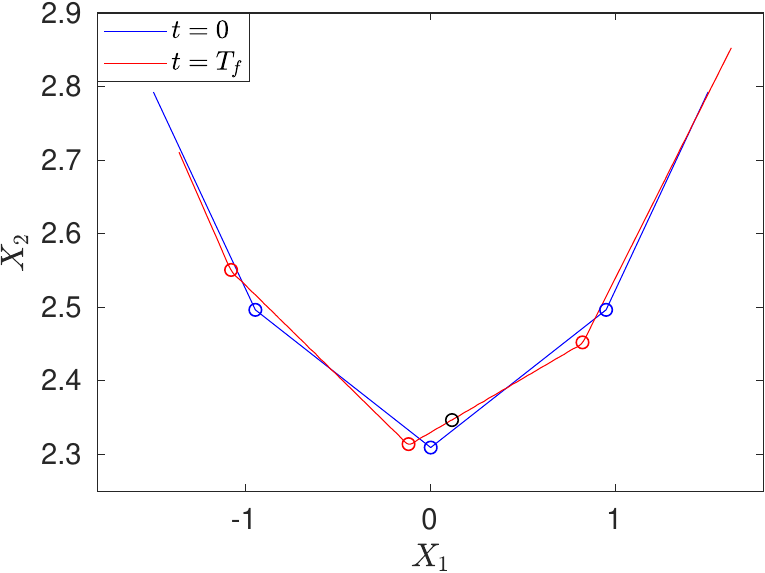}
\includegraphics[width=0.333\textwidth]{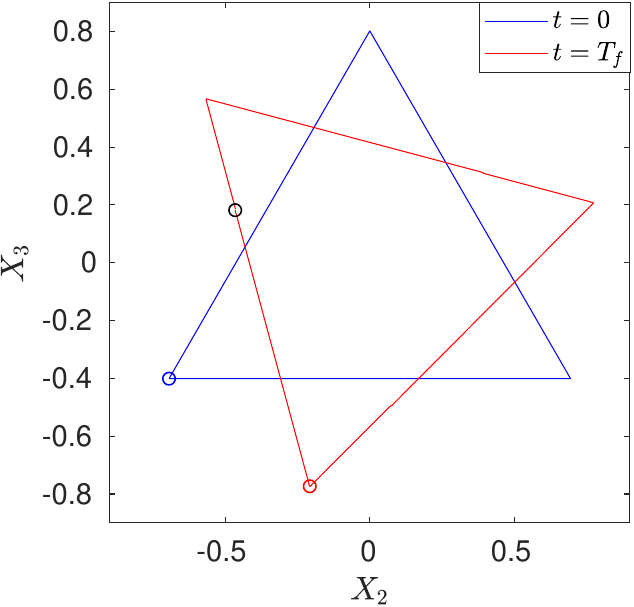}
\caption{The Galilean shift $\spq$ and phase shift: A HHP (left) with $l=0.4$, $\theta_0=\pi/4$, and a CHP (right) for $M=3$, $b=1.2$. The vertices are marked with circles on the curve at time $t=0$ (in blue), $t=T_f$ (in red). Note that at $t=T_f$ the point $s=0$ (in black circle) is shifted by the amount $\spq$ and does not correspond to a corner. Also, the vertices at $t=T_f$ are rotated by a certain amount, i.e., phase shift, with respect to the ones at $t=0$ (in blue).  	
}
\label{fig:Gal-phase-HHP}   
\end{figure} 

Given a rational time $\tpq$, from $\spq$ we know the exact location of corners, as a result, for a curve with $M$ sides at $t=0$, we can compute the angle $\rho_q$ between any two sides of the new helical polygon numerically as 
\begin{equation}
\label{eq:rho_pq_num}
\rho_{pq}^{num,j} = \arccosh(-\T_j \circ_- \T_{j+1}), \ j=0,1,\ldots, Mq-1.
\end{equation}
Here, the value of the piecewise constant tangent vector $\T_j$ has been computed using the mean of the inner points for each side followed by a renormalization. For example, for the CHP, $\T(s)$ with $s\in[0,2\pi/Mq)$, the mean is taken for values corresponding to $s\in [\pi/2Mq,3\pi/2Mq)$ and for the HHP, with $s\in[0,l/q)$, it's taken for $s\in [l/3q,2l/3q)$. Moreover, for the latter, due to the lack of accuracy towards the boundary, we work with only inner $Mq/2$ sides for the comparison; note that a similar approach was used in the planar polygon problem as well, for instance see \cite[Figure 5b]{HozKumarVega2020}. As mentioned previously, the exact value of $\rho_q$ is derived from a conservation law which also implies that the following quantity is preserved: 
\begin{equation}
\label{eq:Ptpq}
P(\tpq) =\prod_m \cosh\left(\frac{\rho_{q}(\tpq)}{2}\right) = \text{constant}, \ m\in \{0,1,\ldots,	\texttt{number of sides}-1 \}. 
\end{equation}
So, to compare $\rho_q$ and $\rho_{pq}^{num,j}$, we have computed the absolute error between $P(\tpq)$ and $P(0)$ for different values of $(p,q)$ for each of the two polygons, and stored them in Table \ref{table:Ptpq-errors}.
 \begin{table}[!htbp]
 	\centering
 	\begin{tabular}{|c|l|l|}
 		\cline{2-3}
 		\multicolumn{1}{c|}{} & \multicolumn{2}{c|}{$|P(t_{pq})-P(0)|$}
 		\\
 		\hline
 		$(p,q)$ & CHP & HHP
 		\\
 		\hline
 		(1,12) & $4.9448\cdot10^{-9}$ & $1.3986\cdot10^{-4}$
 		\\ 		
 		(1,10) & $2.5220\cdot10^{-9}$ & $2.1517\cdot10^{-5}$
 		\\ 		
 		(1,6) & $3.0997\cdot10^{-12}$ & $2.4300\cdot10^{-5}$
 		\\ 		
 		(1,5) & $1.2506\cdot10^{-8}$ & $1.0631\cdot10^{-4}$
 		\\ 		
 		(1,5) & $9.3370\cdot10^{-13}$ & $8.2325\cdot10^{-6}$
 		\\ 		 		
 		(3,10) & $3.6212\cdot10^{-8}$ & $1.1513\cdot10^{-4}$
 		\\		
 		(1,3) & $1.1491\cdot10^{-10}$ & $7.0884\cdot10^{-6}$
 		\\ 				
 		(2,5) & $7.2033\cdot10^{-8}$ & $1.2328\cdot10^{-4}$ 	
 		\\
 		(5,12) & $3.2715\cdot10^{-7}$ & $2.2515\cdot10^{-4}$
 		\\ 				
 		(1,2) & $1.3662\cdot10^{-12}$ & $1.5633\cdot10^{-6}$ 	 		
 		\\	
 		\hline  
 	\end{tabular}~\begin{tabular}{|c|l|l|}
 	\cline{2-3}
 	\multicolumn{1}{c|}{} & \multicolumn{2}{c|}{$|P(t_{pq})-P(0)|$}
 	\\
 	\hline
 	$(p,q)$ & CHP & HHP
 	\\
 	\hline
 	(7,12) & $7.4081\cdot10^{-7}$ & $6.3502\cdot10^{-7}$
 	\\
 	(3,5) & $1.8164\cdot10^{-7}$ & $4.3420\cdot10^{-4}$
 	\\
 	(2,3) & $2.8344\cdot10^{-10}$ & $1.1852\cdot10^{-5}$
 	\\
 	(7,10) & $9.8495\cdot10^{-8}$ & $7.5206\cdot10^{-5}$
 	\\
 	(3,4) & $4.4529\cdot10^{-12}$ & $1.2428\cdot10^{-5}$
 	\\
 	(4,5) & $1.6700\cdot10^{-7}$ & $4.9049\cdot10^{-4}$
 	\\
 	(5,6) & $1.8861\cdot10^{-9}$ & $2.5819\cdot10^{-4}$
 	\\
 	(9,10) & $3.1388\cdot10^{-7}$ & $1.1465\cdot10^{-3}$
 	\\
 	(11,12) & $3.6922\cdot10^{-7}$ & $2.0586\cdot10^{-4}$
 	\\ 
 	(1,1) & $2.0053\cdot10^{-12}$ & $3.6134\cdot10^{-6}$
 	\\
 	\hline  
 \end{tabular}	
\caption{Discrepancy between $P(\tpq)$ and $P(0)$ for a CHP and a HHP for different values of $p$ and $q$, such that $\gcd(p,q)=1$. The parameters for the CHP are $P(0)=1.367631$, for $M=6$, $N/M=7680$, $b=1.2$, and those for the HHP are $P(0)=1.1490\ldots$, $M=48$, $N/M=1920$, $b=0.4$, $l=0.2$. The size of errors clearly implies that the quantity $P(\tpq)$ is preserved during the time evolution in each case.}
\label{table:Ptpq-errors}
 \end{table} 
The parameter values for the CHP are $M=6$, $N/M=7680$, $b=1.2$, $P(0)=1.3676\ldots$, and those for the HHP are $M=48$, $N/M=1920$, $b=0.4$, $l=0.2$, $P(0)=1.1490\ldots$. The size of errors clearly implies the agreement between the two approaches; one can also note that the error values are much smaller in the case of a CHP. This follows from the fact that the corresponding trigonometric functions are approximated with the pseudo-spectral method and for a HHP, a finite difference scheme is used for the hyperbolic functions. Moreover, because of the efficiency of the \texttt{fft} algorithm in MATLAB and symmetries of the CHP, working with a large value of $N/M$ is possible unlike for the HHP case (see also the discussion in \cite[Section 3.1]{HozKumarVega2020}).  

Numerical simulations show that apart from the formation of new sides, the polygonal curve also moves in the vertical direction. To measure that, we compute the center of mass by calculating the mean of all the values of $\X(s_k,t^{(n)})$ and call it $\mathbf h(t^{(n)})$. For a CHP, the first component of $\mathbf h(t^{(n)})$ behaves as a linear function of $t$ with a negative slope which corresponds to the speed $c_M$ of the center of mass along the negative $x$-axis; however, the second and third components are zero. On the other hand, for a HHP, bearing in mind \cite[(40)]{HozKumarVega2020}, we continue to work with the inner points $N/2$ nodes of the curve $\X$ and observe that the first two components of corresponding $\mathbf h(t^{(n)})$ are nonzero while the third component grows linearly in time $t$ with a rate $c_l$. 

We compute $c_M$, $c_l$ numerically, i.e., $c_M^{num}$, $c_l^{num}$, respectively, and compare them with their corresponding exact values  
\begin{equation}
\label{eq:c_M-circ-hel}
c_M = \frac{2\ln(\cosh(\rho_0/2))}{(\pi/M)\tanh(\pi/M)}, \ c_l = \frac{2\ln(\cosh(\rho_0/2))}{(l/2)\tanh(l/2)},
\end{equation}
which can be obtained using the techniques described in \cite[Section 4]{HozVega2018}. For numerical simulations, we have taken $M=6,7,\ldots,15$, discretization $N/M=480,960,\ldots,7680$, $b=1.2$ for the CHP, and for the HHP $M=48$, $b=0.4$, $l=0.08, 0.1, \ldots, 0.2$, $N/M=240, 480, \ldots, 1920$. The relative errors in each case are displayed on the left-hand side of Figures \ref{fig:cM-error} and \ref{fig:cl-error}, while the right-hand side shows the plot of $c_M$ for different $M$ values which approaches to $b^2-1=0.44$ as $M$ gets larger, and $c_l$ tending to $b^2+1=1.16$ as $l$ becomes smaller, respectively (plotted in dash-dotted red line).   
 
In both cases, the numerical experiments clearly show that as the discretization is doubled the relative errors decrease by a factor slightly less than two, thus, a first-order convergence. 
	
\begin{figure}[htbp!] 
\centering
\includegraphics[width=0.45\textwidth]{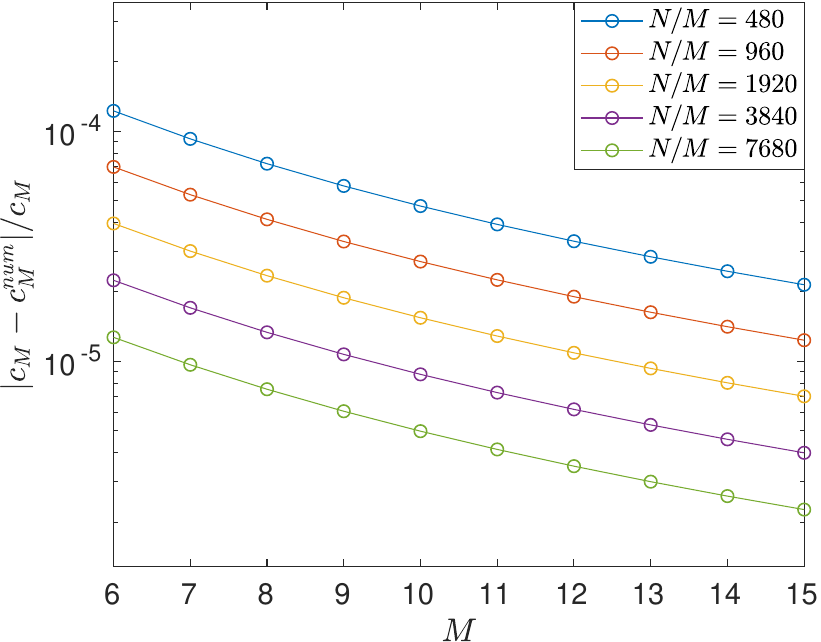}
\includegraphics[width=0.45\textwidth]{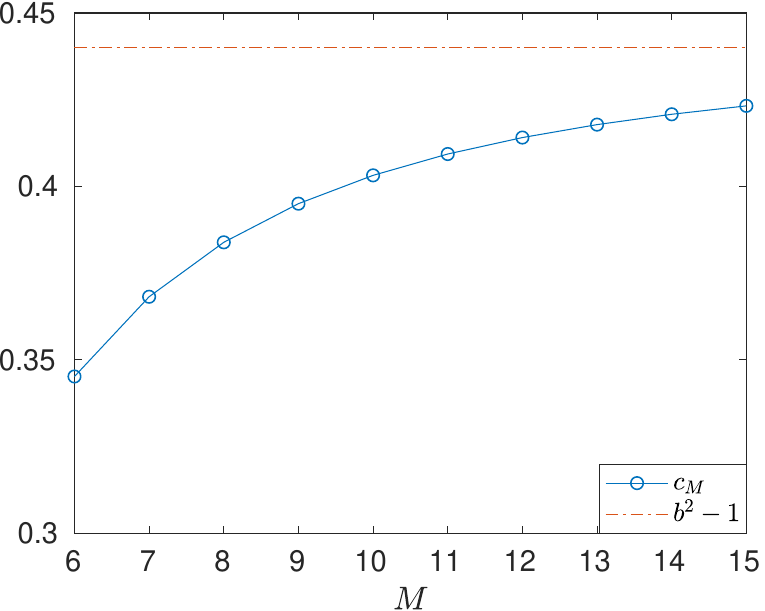}
\caption[$|c_M-c_M^{num}|/c_M$ for $b=1.2$ and different $M$ and discretizations.]{Left: $|c_M-c_M^{num}|/c_M$, computed for $b=1.2$, and different values of $M$ and $N/M$. The error clearly decreases, as $N/M$ increases, showing the convergence to the theoretical value. Right: $c_M$, computed using \eqref{eq:c_M-circ-hel}, for different values of $M$. We have also plotted in dash-dotted line the limiting value $b^2-1=0.44$, to which we conjecture $c_M$ to tend, as $M\to\infty$.}
\label{fig:cM-error}   
\end{figure}
\begin{figure}[htbp!] 
\centering
\includegraphics[width=0.48\textwidth,valign=b]{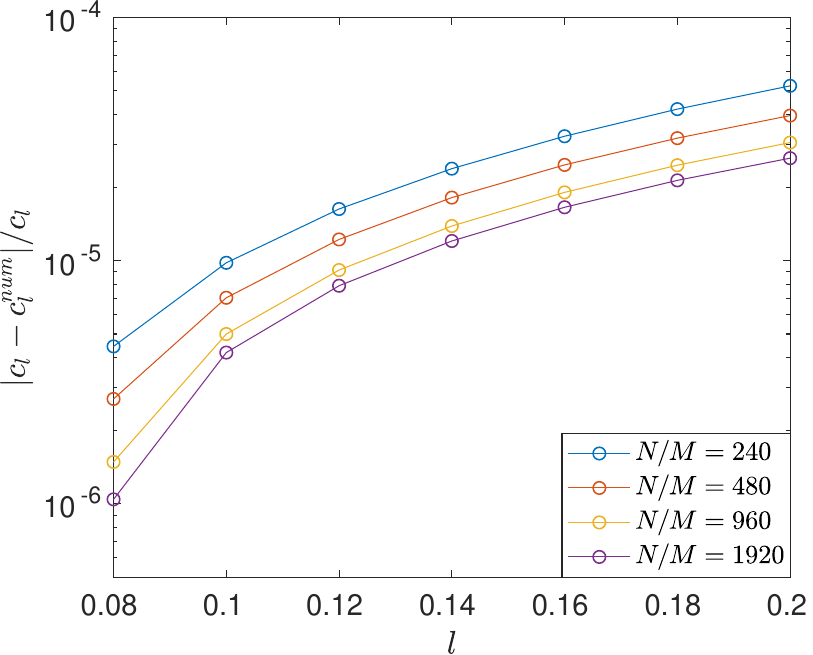}
\includegraphics[width=0.455\textwidth,valign=b]{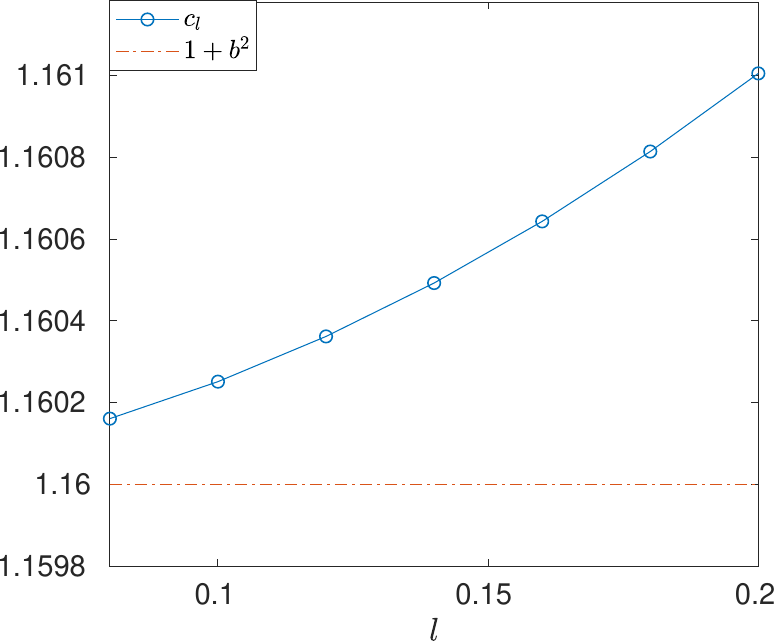}
\caption[$|c_l-c_l^{num}|/c_l$ for $b=0.4$ and different $l$ and discretizations.]{Left: $|c_l-c_l^{num}|/c_l$, computed for $b=0.4$, $M=48$, and different values of $l$ and $N/M$. The errors follow a similar pattern as in Figure \ref{fig:cM-error}. Right: $c_l$, computed using \eqref{eq:c_M-circ-hel}, for different values of $l$. The limiting value $1+b^2=1.16$, to which we conjecture $c_l$ to tend, as $l\to0$ is plotted in dash-dotted line.}
\label{fig:cl-error}   
\end{figure}
\section{Trajectory of $\X(0,t)$}
\label{sec:Traj-X0t}
In the evolution of \eqref{eq:SMP-hyp}--\eqref{eq:VFE-hyp} for regular polygons, trajectory of a single point has revealed the most striking phenomenon of multifractality \cite{HozVega2014,HozKumarVega2019,HozKumarVega2020}. When $\X(s,0)$ is a planar $l$-polygon, due to its symmetries, the evolution of one point, i.e., $\X(0,t)$ which through numerical simulation was claimed to be multifractal, lies in the YZ-plane. However, it's no more planar when the torsion is introduced in the initial data. In the following lines, we explore three different cases determined by the choice of the parameter $b$. 
\subsection{When $b\in(1,\infty)$}
Corresponding to a CHP, for a given values of $M$ and $b$ $(\text{or}, \theta_0)$, through numerical simulations, a non-planar multifractal structure is observed in $\X(0,t)$. To understand it further, for a time-period larger than $T_f$ (unlike in the zero-torsion case), for parameters $M=6$, $\theta_0 = 2\pi/5$, $t\in[0,10\pi/3]$, we note that $\X(0,t)$ consists of a fractal structure which is repeated periodically with a certain counterclockwise rotation and a translation along the $x$-axis. This behaviour, similar to the one in the Euclidean helical polygons, is shown in Figure \ref{fig:X0t-M6-hyp} where $\X(0,t)$ (in black) is plotted along with $\X(s,0)$ (in blue) and $\X(s,10\pi/3)$ (in red) \cite{HozKumarVega2019}. Then, by following the discussion in \cite[Section 4.1]{HozKumarVega2019}, it is clear that when $\theta_0>0$, the occurrence of a corner in the curve $\X(s,t)$ and thus, in $\X(0,t)$ is determined by the Galilean shift $\spq$. As a result, by choosing $\theta_0 = c \pi / d$, $\gcd(c,d)=1$, after $d$ or $d/2$ time periods, $s_{pq}$ becomes equal to the side-length of the helical polygon and the point $s=0$ and thus, $\X(0,t)$ has a corner. More precisely, if we define 
\begin{equation}
\label{eq:Tfcd}
T_f^{c,d} = 
\begin{cases}
(d/2)\, T_f \ &\text{if $c\cdot d$ odd}, \\
d \, T_f \ &\text{if $c\cdot d$ even},
\end{cases}
\end{equation}
then at every integer multiples of $T_f^{c,d}$, $\X(0,t)$ will have a corner with a repetition of multifractal structure, and at the rational multiples of $T_f^{c,d}$, corners of different scale would appear. 
\begin{figure}[!htbp]\centering
\includegraphics[width=0.6\textwidth, clip=true]{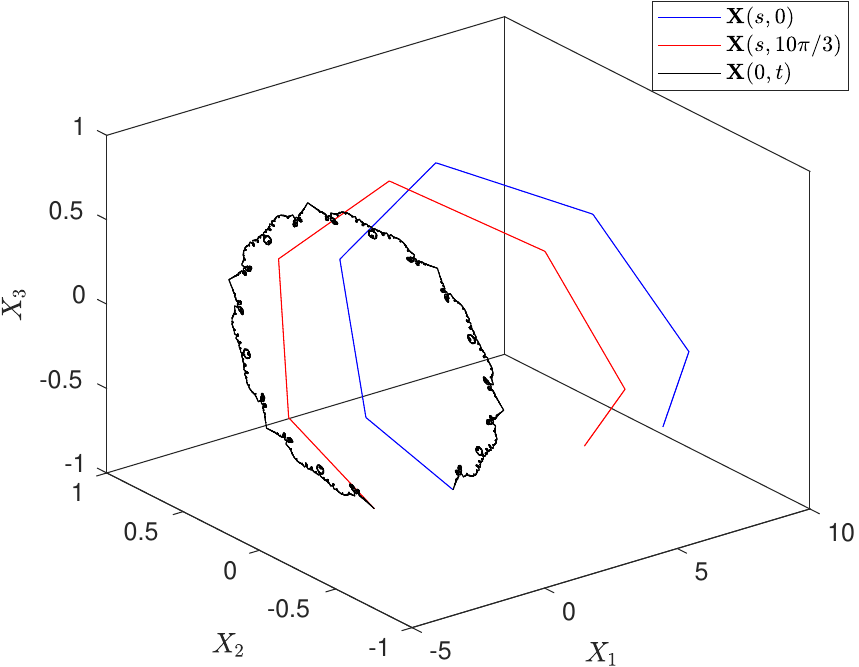}
\caption{$\X(s,0)$ (in blue) for a CHP, with parameters $M=6$, $\theta_0=2\pi/5$, and its evolution at time $t =  10\pi/3$ (in red), together with the curve described by $\X(0,t)$, for $t\in[0,10\pi/3]$ (in black). $\X(0,t)$ has a conspicuous fractal structure which repeats periodically, with some rotation and a vertical movement.}\label{fig:X0t-M6-hyp}
\end{figure}
 
In order to further understand $\X(0,t)$, for each of its components, we compute its \textit{fingerprint}, i.e., a plot of the Fourier coefficients and wavenumbers which was found very useful in understanding the multifractal patterns in \cite{HozKumarVega2019}. We define 
\begin{equation}
\label{eq:Rt-nut-def-CHP}
\begin{cases}
z_{2,3}(t) = X_2(0,t) + i X_3(0,t) = R(t) e^{i\nu(t)}, \\
\tilde X_1(t) = X_1(0,t) + c_M t, \ t\in[0,T_f^{c,d}],
\end{cases}
\end{equation}    	 
with $R(t) = \sqrt{X_2(0,t)^2 + X_3(0,t)^2}$ and $\nu(t) = \arctan(X_3(0,t)/X_2(0,t))$. Since both $R(t)$ and $\tilde X_1(t)$ are periodic, we express them in their Fourier series expansion and plot the product of corresponding Fourier coefficients and wavenumbers against wavenumbers. We refer to this as a fingerprint plot. In Figure \ref{fig:Rt-M6-hyp}, for parameters $M=6$, $\theta_0=2\pi/5$, $N/M=2^{11}$, on the left-hand side we have $R(t)$, its fingerprint plot on the right-hand side where $a_{n,M}$, $n=1,2,\ldots,2000$ are its Fourier coefficients that are approximated using \texttt{fft} algorithm in MATLAB. On the left-hand side of Figure \ref{fig:nu_t-M6-hyp}, we plot $\nu(t)$, which besides a conspicuous fractal structure exhibits a linear behaviour in $t$ with a slope equal to $-1.1715\ldots$. If we remove the linear part in $\nu(t)$, then the resulting curve, let's call it $\tilde \nu(t)$, is $T_{f}^{c,d}$-periodic and we show its fingerprint plot on the right-hand side of Figure \ref{fig:nu_t-M6-hyp}.
\begin{figure}[!htbp]\centering
\includegraphics[width=0.472\textwidth, clip=true]{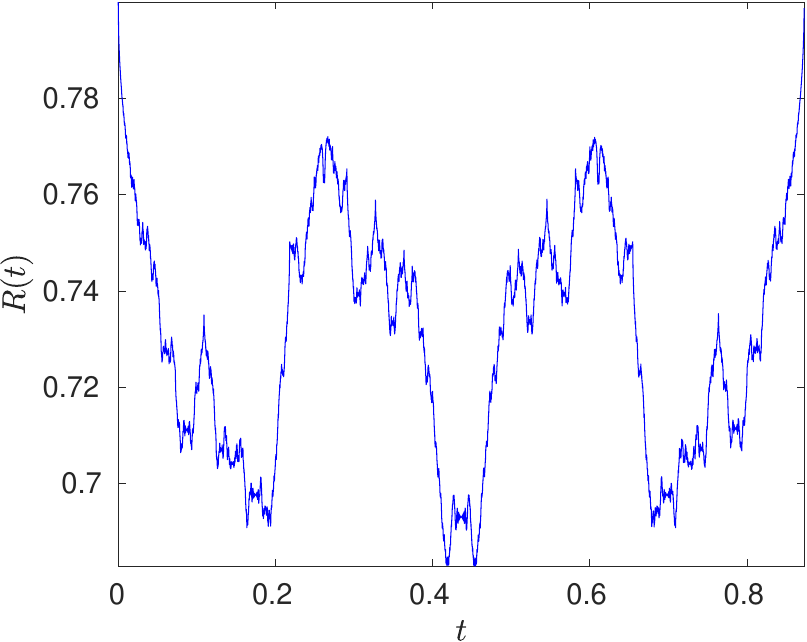}
\includegraphics[width=0.497\textwidth, clip=true]{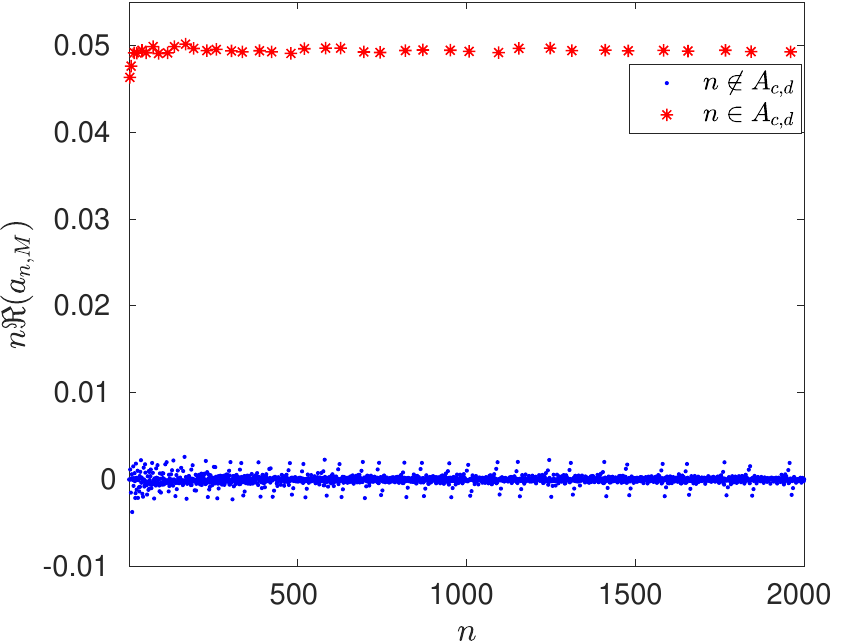}
\caption{For the CHP with parameters $M=6$, $\theta_0=2\pi/5$, $N/M=2^{11}$; left: $R(t)$, for $t\in[0,T_f^{c,d}]$. Right: the fingerprint plot of $R(t)$, where $a_{n,M}$ are its Fourier coefficients and the dominating points (red starred) are given by \eqref{eq:Acd}. }\label{fig:Rt-M6-hyp}
\end{figure}
\begin{figure}[!htbp]\centering
	\includegraphics[width=0.445\textwidth, clip=true]{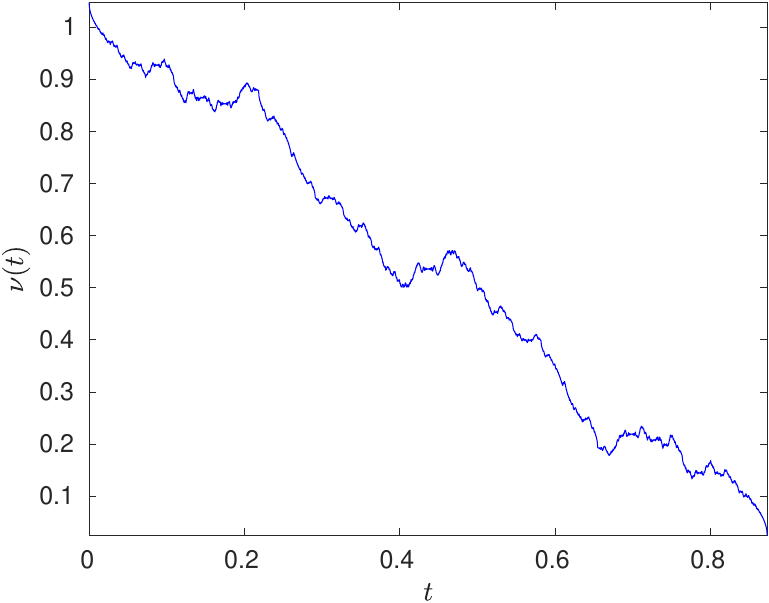}
	\includegraphics[width=0.485\textwidth, clip=true]{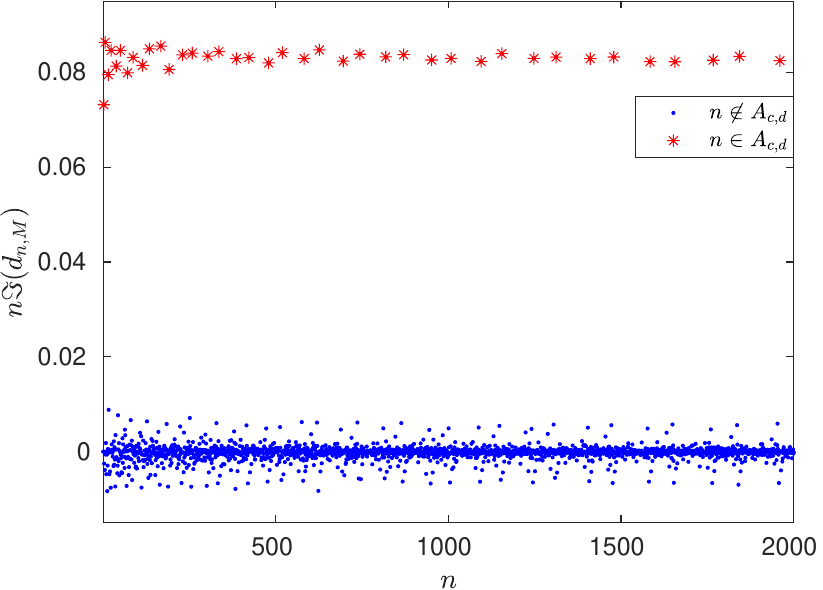}
	\caption{For the CHP with parameters $M=6$, $\theta_0=2\pi/5$, $N/M=2^{11}$; left: $\nu(t)$, for $t\in[0,T_f^{c,d}]$, which, besides a conspicuous fractal structure, also varies linearly in $t$ with a slope equal to $-1.1715\ldots$. Right: the fingerprint plot of $\tilde \nu(t)$, i.e., $\nu(t)$ without the linear part, where $d_{n,M}$ are its Fourier coefficients and the dominating points (red starred) are given by \eqref{eq:Acd}. }\label{fig:nu_t-M6-hyp}
\end{figure}

We have performed several numerical experiments to understand $R(t)$, $\tilde X_1(t)$, $\nu(t)$, and the fingerprint plots reveal that when the initial torsion angle $\theta_0=c \pi /d$, $\gcd(c,d)=1$, the dominating points of the fingerprint plot belong to the set
\begin{equation}
\label{eq:Acd}
A_{c,d}=
\begin{cases}
\{n(nd+c)/2 \, | \, n\in \mathbb{Z}\} \cap \mathbb{N} \ &\text{if $c\cdot d$ odd}, \\
\{n(nd+c) \, | \, n\in \mathbb{Z}\} \cap \mathbb{N} \ &\text{if $c\cdot d$ even}.
\end{cases}
\end{equation}

Next, the connection between $\X(0,t)$ for regular polygons and Riemann's non-differentiable function, prompts us to compare $\tilde X_1(t)$ and $\tilde \nu(t)$, with the imaginary part of 
\begin{equation}
\label{eq:phi_cd}
\phi_{c,d}(t) = \sum_{k\in A_{c,d}} \frac{e^{2\pi i kt}}{k}, \ 
t\in
\begin{cases}
[0,1/2] \ &\text{if $c\cdot d$ odd}, \\
[0,1] \ &\text{if $c\cdot d$ even}.
\end{cases}
\end{equation}
For the same parameter values, the left side of Figure \ref{fig:XM1-M6} shows the curves $-\tilde X_1(t)$, $\Im(\phi_{c,d}(t))$ with $c=2$, $d=5$, and one can note the similarities between the two figures. On the right-hand side, we have the fingerprint of the scaled $\tilde X_1(t)$ with $b_{n,M}$ as the approximation of its Fourier coefficients. The dominating points shown in red color belong to the set $\Acd$ and can be seen distributed around the value $1/2$ while the rest around the value 0. A similar behaviour is observed for $\tilde{\nu}(t)$ as well. 
\begin{figure}
\centering
\begin{minipage}{0.472\linewidth}
\begin{subfigure}{\linewidth}
\includegraphics[width=\linewidth]{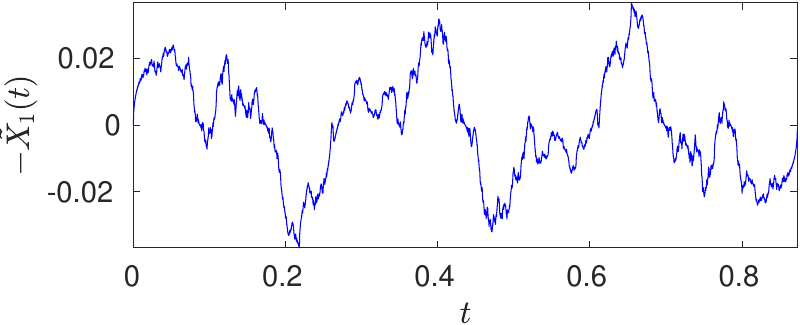}
\end{subfigure}	
\begin{subfigure}{\textwidth}
\includegraphics[width=\linewidth]{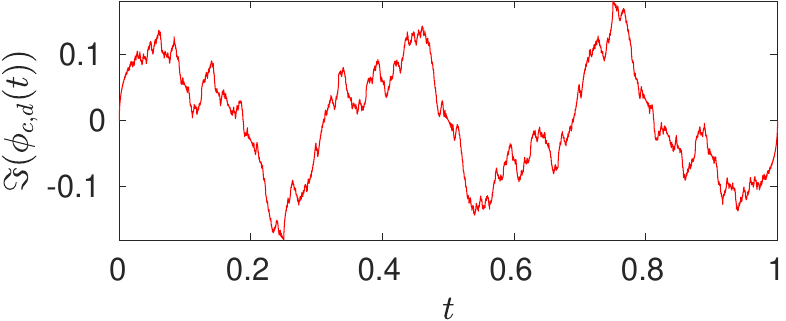}
\end{subfigure}
\end{minipage}
\hfil
\begin{minipage}{0.52\linewidth}
\begin{subfigure}{\linewidth} 	
\includegraphics[width=\linewidth]{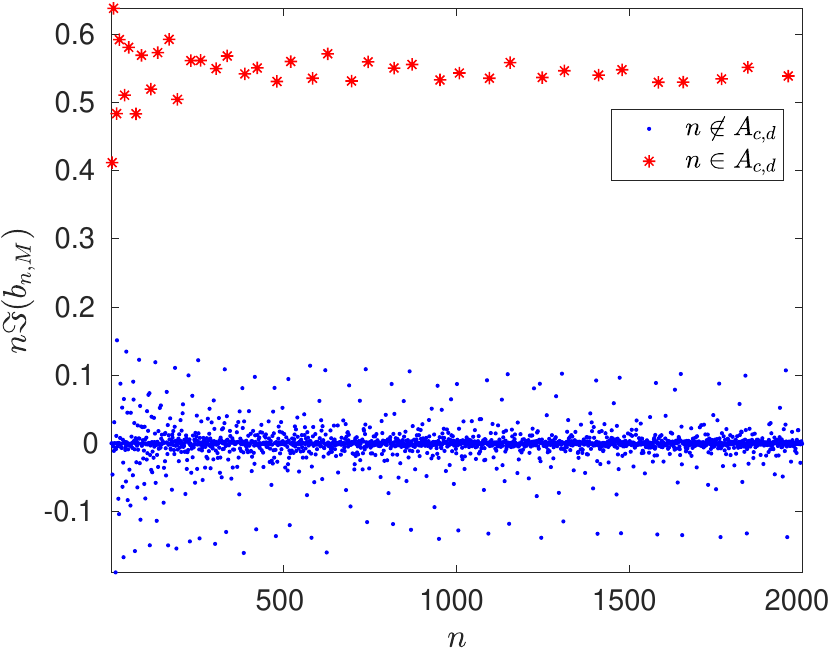}
\end{subfigure}
\end{minipage}
\caption{For a CHP with parameters $M=6, \theta_0 = 2\pi/5,$ $N/M=2^{11}$; left: on the top is $-\tilde{X}_1(t)$, $t\in[0,\tcd]$ and $\Im({\phi_{c,d}(t)})$, $t\in[0,1]$ at the bottom.  Right: the fingerprint plot of the scaled $\tilde X_1$, where $b_{n,M}$ are its Fourier coefficients, and the dominating points (red starred) are given by \eqref{eq:Acd}.}
\label{fig:XM1-M6}
\end{figure}
\subsection{When $b\in(0,\infty)$} The numerical simulations for a HHP that is characterized with parameters $l$, $b$ and $M$, show that the multifractal behaviour in $\X(0,t)$ repeats periodically with a certain counterclockwise rotation and a translation along the $z$-axis. As before, thanks to the Galilean shift $s_{pq} = l \theta_0 p / (\pi q)$, by choosing $\theta_0= \pi c /d $, $\gcd(c,d)=1$, the periodicity in space and time can be recovered. In other words, for $\tcd$ as in \eqref{eq:Tfcd} with $T_f=l^2/2\pi$, if we define 
\begin{equation}
\label{eq:Rt-nut-def-HHP}
R(t) = \sqrt{-X_1(0,t)^2+X_2(0,t)^2}, \quad \nu(t) = \arctanh\left(\frac{X_1(0,t)}{X_2(0,t)}\right),  \quad
\tilde X_3(t) = X_3(0,t)-c_l t,
\end{equation}
then for parameters $M=96$, $l=0.2$, $\theta_0=\pi/4$, $N/M=2^{10}$, $t\in[0,\tcd]$, we notice that $R(t)$ plotted on the left-hand side of Figure \ref{fig:M96-Rt-nu} and $\tilde X_3(t)$ on the top left-hand side of Figure \ref{fig:M96-XM3} are $\tcd$-periodic. Furthermore, as in the case of a CHP, $\nu(t)$ has a multifractal structure with a linear behaviour in $t$ with a slope equal to $3.8068\ldots$, plotted in the center of Figure \ref{fig:M96-Rt-nu}, and on the right-hand side is $\tilde \nu(t)$, i.e., $\nu(t)$ without the linear part.
The fingerprint plots of these periodic curves show that the dominating points belong to the set $\Acd$ in \eqref{eq:Acd}; consequently, up to a scaling, $\tilde{\nu}(t)$ and $\tilde X_3(t)$ can be compared with the imaginary part of $\phi_{c,d}$ in \eqref{eq:phi_cd}, see Figure \ref{fig:M96-XM3}.
\begin{figure}[!htbp]\centering
	\includegraphics[width=0.32\textwidth, clip=true]{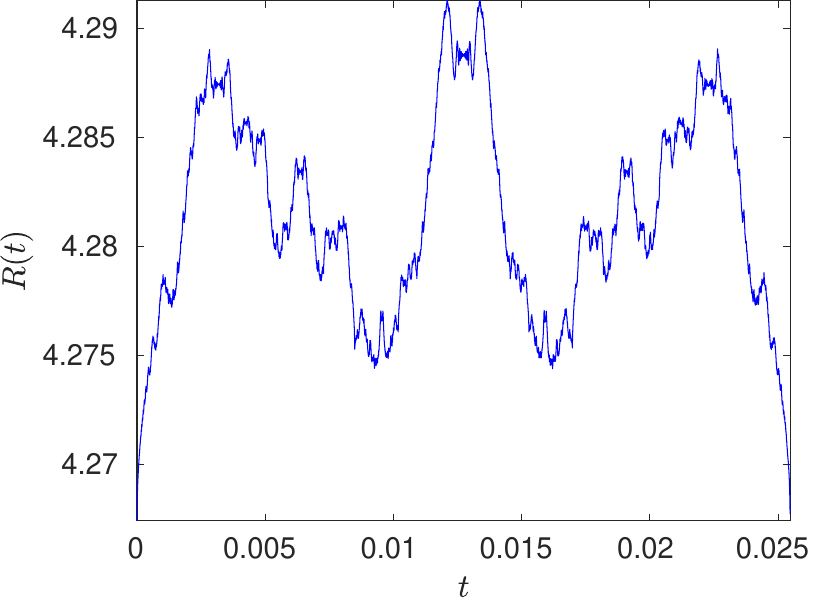}
	\includegraphics[width=0.32\textwidth, clip=true]{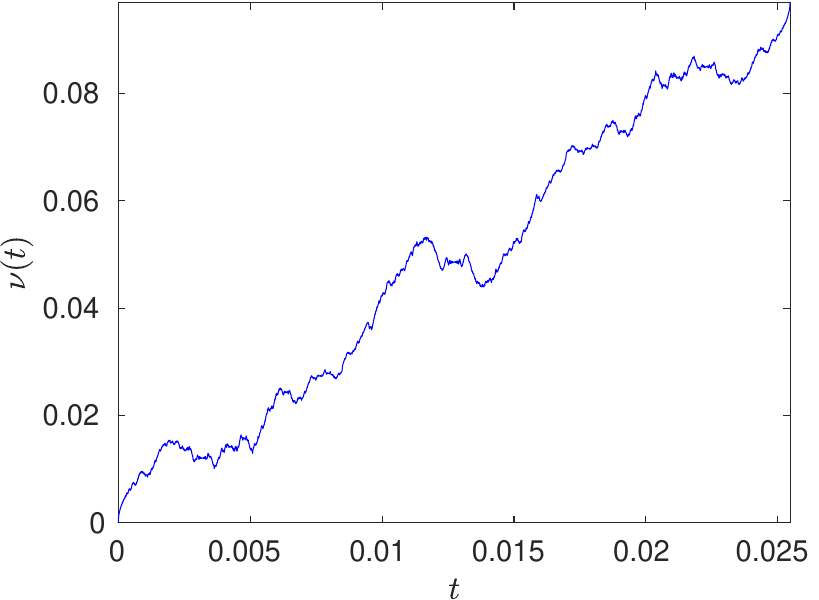}
	\includegraphics[width=0.31\textwidth, clip=true]{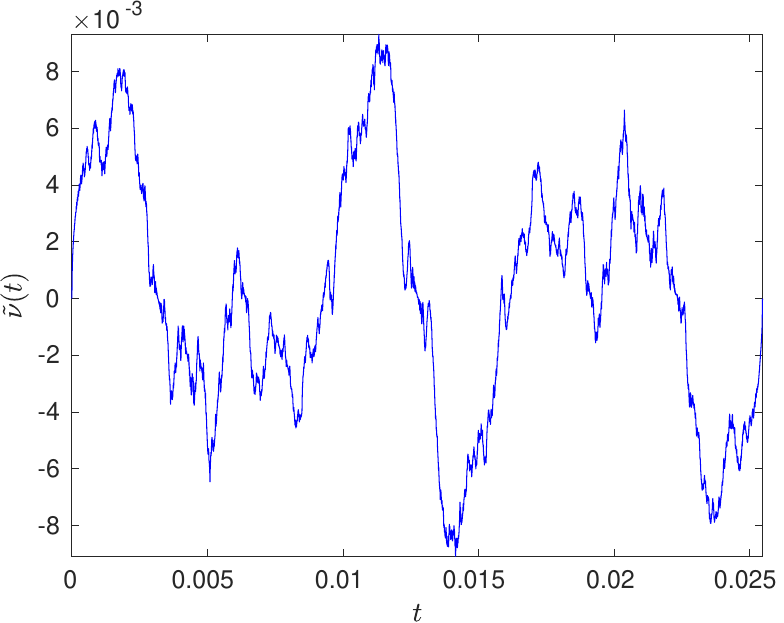}
	\caption[]{For a HHP with parameters $M=96$, $\theta_0=\pi/4$, $l=0.2$, $N/M=2^{10}$; left: $R(t)$, $t\in [0,T_f^{c,d}]$. Center: $\nu(t)$ as defined in \eqref{eq:Rt-nut-def-HHP}, which, besides a conspicuous fractal structure, also varies linearly in $t$ with a rate equal to $3.8068\ldots$. Right: The periodic curve $\tilde{\nu}(t) = \nu(t)-3.8068 \ t$, which, up to a scaling, can be compared with the imaginary part of $\phi_{c,d}$ in Figure \ref{fig:M96-XM3}.}
	\label{fig:M96-Rt-nu}
\end{figure}			
\begin{figure}
	\centering
	\begin{minipage}{0.49\linewidth}
		\begin{subfigure}{\linewidth}
			\includegraphics[width=\linewidth]{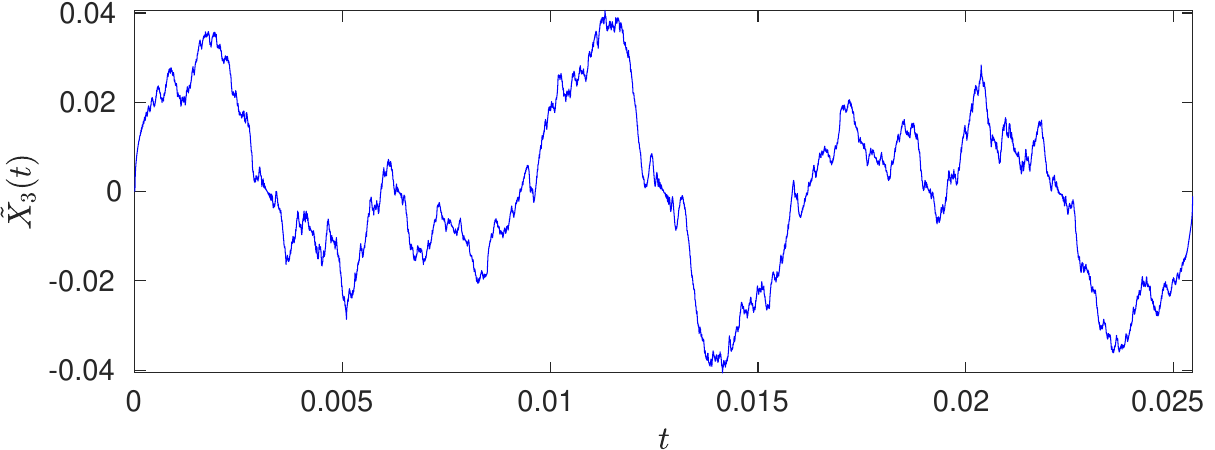}
		\end{subfigure}	
		\begin{subfigure}{\textwidth}
			\includegraphics[width=\linewidth]{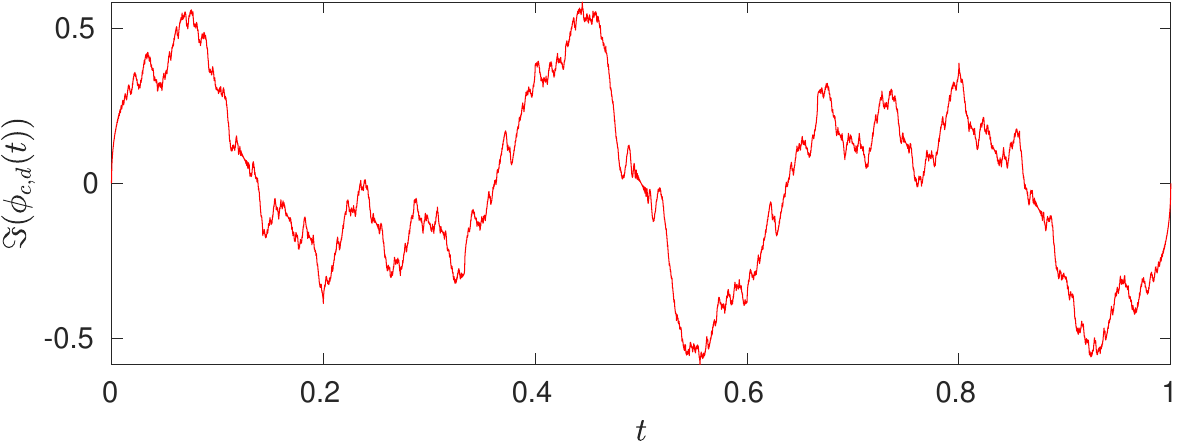}
		\end{subfigure}
	\end{minipage}
	\hfil
	\begin{minipage}{0.50\linewidth}
	\begin{subfigure}{\linewidth}
		\includegraphics[width=\linewidth]{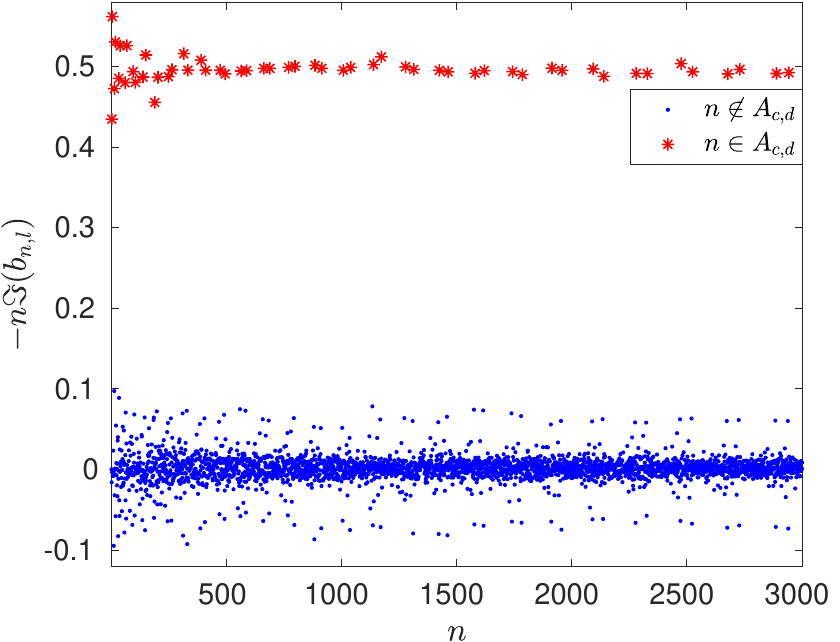}
	\end{subfigure}
\end{minipage}
	\caption{For a HHP with parameters $M=96, \theta_0 = \pi/4,$ $N/M=2^{10}$, $l=0.2$; left: On the top is $\tilde{X}_3(t)$, $t\in[0,\tcd]$ and bottom is $\Im{(\phi_{c,d}(t))}$, $t\in[0,1]$.  Right: Approximation of $n\Im(b_{n,l})$, for $n=1,2,\ldots, 2000$, where $b_{n,l}$ are the Fourier coefficients of the scaled $\tilde X_3(t)$. The dominating points (red starred) are given by \eqref{eq:Acd}}
	\label{fig:M96-XM3}
\end{figure}
\subsection{When $M\to\infty$, or $l\to0$}
In case of a CHP, for a fixed value of $b$, when the value of $M$ is increased, the initial curve $\X$ tends to a circular helix; see Figure \ref{fig:M20-X0t-hyp} where $M=20$, $\theta_0=\pi/9$, and $\X(0,t)$ (in black), $t\in[0,9\pi/10]$, is plotted with $\X(s,0)$ (in blue) and $\X(s,9\pi/10)$ (in red). Moreover, the fingerprint plot of the scaled $\tilde X_1(t)$ on the left-hand side of Figure \ref{fig:M20-XM1-nu} and that of scaled $\tilde{\nu}(t)$ on the right-hand side show that the dominating points belonging to the set $\Acd$ are more orderly distributed around the value $1/4$ and the rest around the value 0 when compared with Figure \ref{fig:XM1-M6}. Here, the slope of linearity in $\nu(t)$ is equal to $-1.1100$, whose absolute value is very close to $b=1.1133$.

A similar behaviour is observed for a HHP when $l$ is very small, for instance, see the left-hand side of Figure \ref{fig:M96-XM3-nu-l005} for $l=0.05$, $\theta_0=\pi/16$, $M=96$, $N/M=2^{10}$, where the dominating points are distributed around the value $1/2$. Here, the slope of linearity in $\nu(t)$ is equal to $3.9212\ldots$ which is very close to the value of $b=3.9405\ldots$. 

Thus, after checking the fingerprint plots for several parameter values in each case, the consistency in the numerical results allow us to conjecture that 
\begin{equation}
\label{eq:bnM_large}
\lim |n \, b_{n}|=
\begin{cases}
1/4 \ &\text{if $c\cdot d$ odd}, \\
1/2 \ &\text{if $c\cdot d$ even}, \\
0 \ &\text{otherwise},
\end{cases}
\end{equation}
where the Fourier coefficients $b_n\equiv b_{n,M}$, for a CHP and limit is as $M\to\infty$, and $b_n\equiv b_{n,l}$, for a HHP and limit is as $l\to0$. The conjecture \eqref{eq:bnM_large} holds true also for $\nu$, where using the same notations, instead of $b_{n}$, we will have $d_{n}$, i.e., the Fourier coefficients of $\tilde \nu$, thus, indicating that $\tilde \nu$ converges to $\Im{(\phi_{c,d})}$. 	

\begin{figure}[!htbp]\centering
	\includegraphics[width=0.7\textwidth, clip=true]{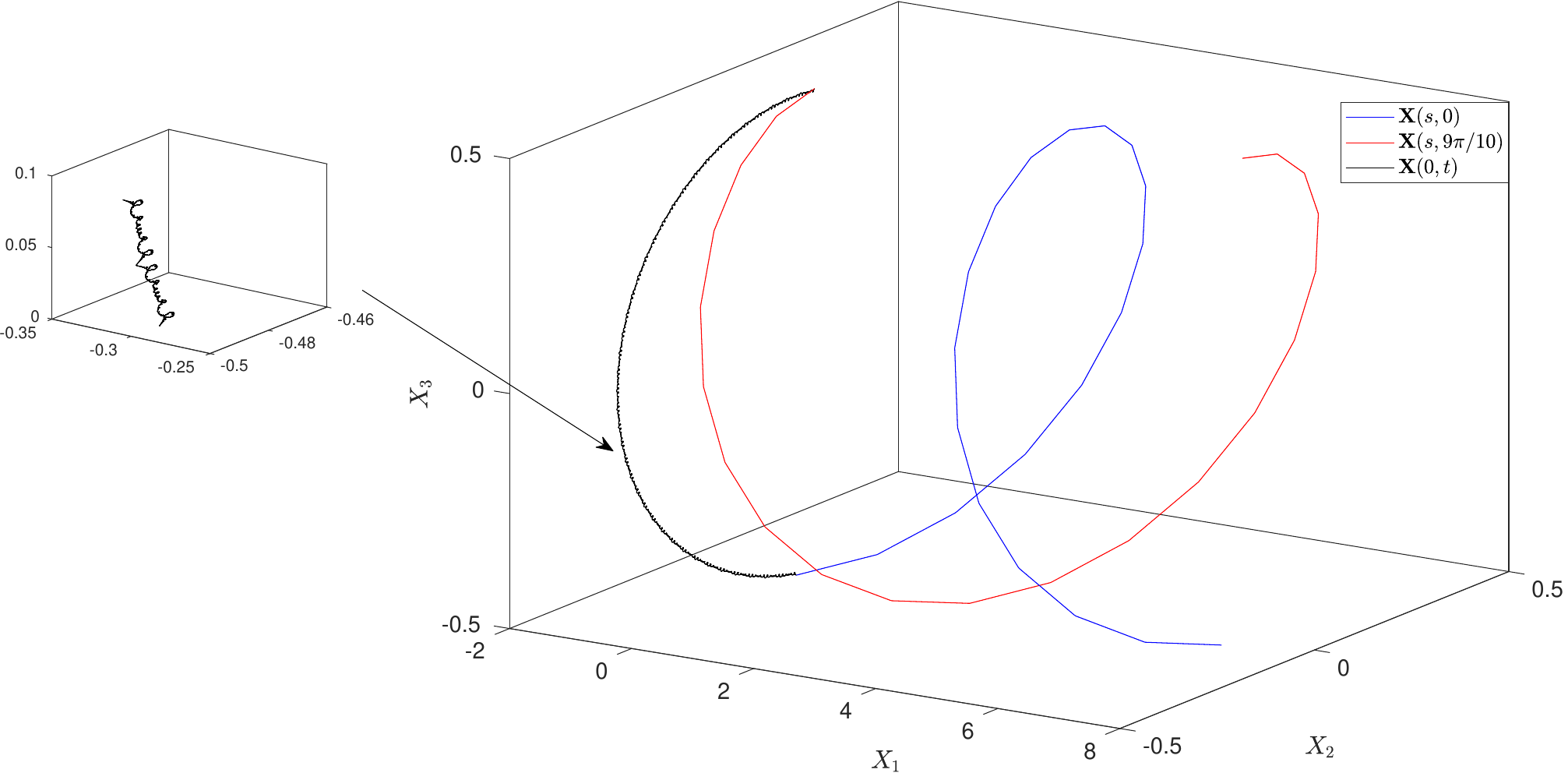}
	\caption{$\X(s,0)$ (in blue) for a CHP with parameters $M=20$, $\theta_0=\pi/9$, and its evolution at time $t  = 9\pi/10$ (in red), together with the curve described by $\X(0,t)$, for $t\in[0,9\pi/10]$ (in black). $\X(0,t)$ has a fractal helical structure, as shown in the zoomed part.}
	\label{fig:M20-X0t-hyp}
\end{figure}
\begin{figure}[!htbp]\centering
	\includegraphics[width=0.5\textwidth, clip=true]{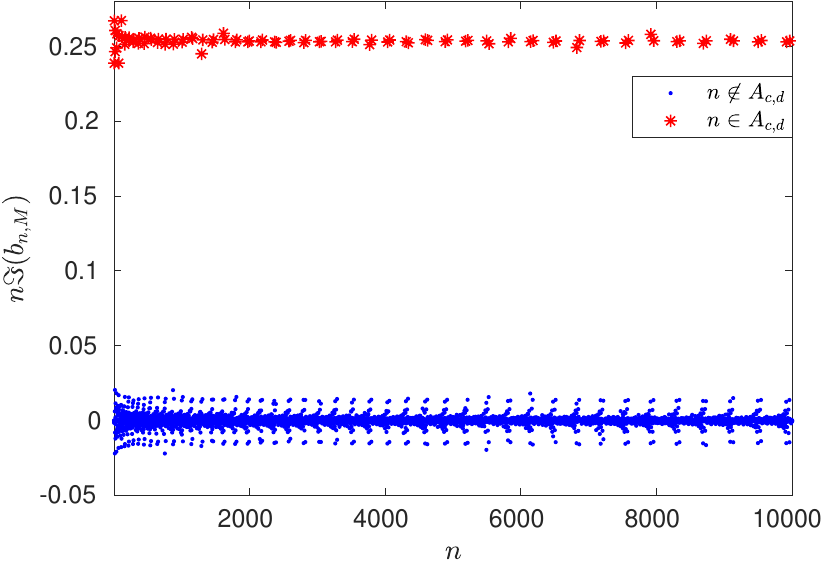}
	\includegraphics[width=0.485\textwidth, clip=true]{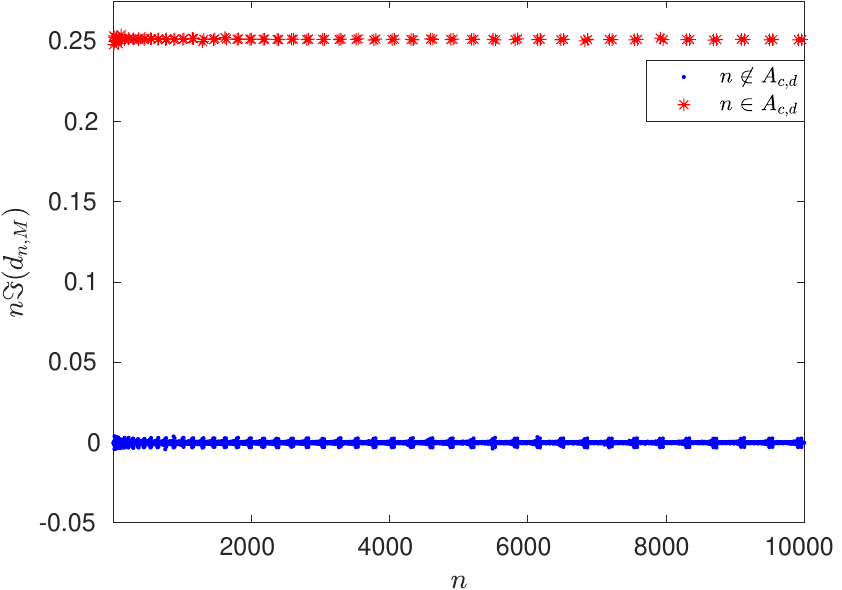}
	\caption{For a CHP with $M=20$, $\theta_0=\pi/9$, $c=1$, $d=9$, $b=1.1133\ldots$; left: the fingerprint of the scaled $\tilde X_{1}(t)$ for , $t\in [0,T_f^{c,d}]$. Right: The fingerprint plot of the scaled $\tilde \nu(t)$. The convergence of the dominating points to the value $1/4$ is visible in both plots; however, it is clearly stronger for $\tilde{\nu}(t)$ as the remaining points (in blue) are much closer to zero.}
	\label{fig:M20-XM1-nu}
\end{figure}			
\begin{figure}[!htbp]\centering
	\includegraphics[width=0.47\textwidth,  valign=t, clip=true]{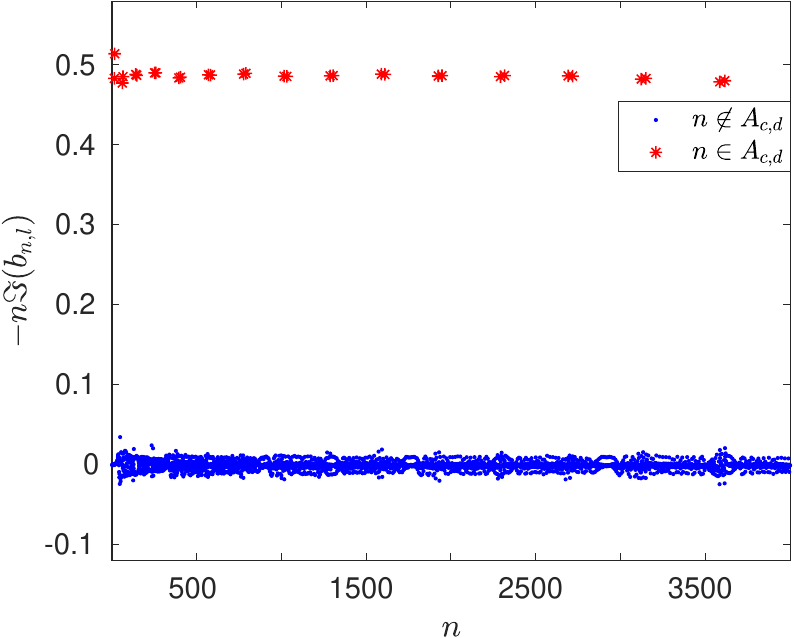}
	\includegraphics[width=0.515\textwidth, valign=t, clip=true]{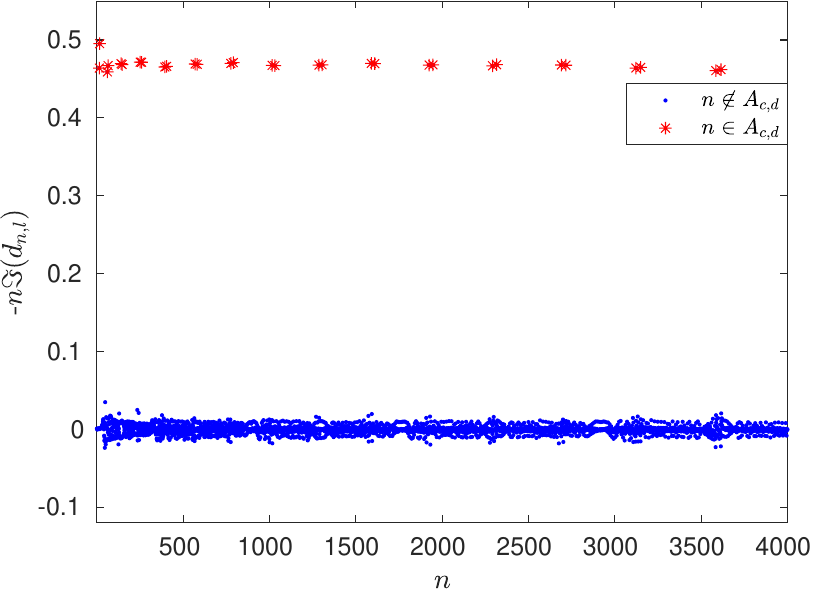}
	\caption{For a HHP with $M=96, \theta_0 = \pi/16,$ $N/M=2^{10}$, $l=0.05$; left: the fingerprint plot of the scaled $\tilde X_3(t)$. Right: The fingerprint plot of the scaled $\tilde \nu(t)$. The convergence of the dominating points (red starred) to the value $1/2$ is clearly visible and the remaining ones are much closer to the value $0$ when compared with Figure \ref{fig:M96-XM3}.}
	\label{fig:M96-XM3-nu-l005}
\end{figure}

Note that $\nu(t)$ can be associated with the phase shift which unlike Galilean shift is rather more complicated and we have analysed it by dividing it into linear and periodic parts. While the linear part seems to converge to $b\ t$, the periodic part to a multifractal $\Im(\phi_{c,d})$, up to a scaling. At this point, let's recall the time evolution of the smooth helical curves discussed in \cite[Section 1.4]{KumarPhd}, where the parameter $b$ also corresponds to the angular speed. Hence, for the CHP and HHP, the slope of the linear part in $\nu(t)$, can be associated to the rate with which the polygonal curve rotates in the YZ-plane and XY-plane, respectively, which in turn converges to $b$ with $M$, or $l$. A precise expression for the slope and the periodic curve $\tilde \nu(t)$ for any $M$, or $l$ is certainly desired; however, appears to be very challenging, which we plan to address in an upcoming work.


\subsection{When $b\to1^+$}
As mentioned before, for a fixed $M$, when the torsion angle $\theta_0 \to 2\pi/M$, $b \to 1^+$, and the initial curve $\X(s,0)$ approaches to a straight line. Thus, by taking $b=1+10^{-5}$ for the numerical simulations, we note that for $t\in[0,2\pi]$ both space and time periodicity are recovered as commented in Section \ref{sec:Pb-f-CHP}. Moreover, as $b\to1^+$, from \eqref{eq:c_M-circ-hel}, $c_M \to 0$, consequently, the movement of $\X(0,t)$ in the YZ-plane is much larger than along the $x$-axis. Therefore, after removing the vertical height from the first component and performing a stereographic projection of the resulting curve, we define 
\begin{equation}
\label{eq:zt}
z(t) = - \frac{X_2(0,t)}{1+\tilde X_1(t)}+ i \frac{X_3(0,t)}{1+\tilde X_1(t)}, \ t \in[0,2\pi],
\end{equation}
which is almost $2\pi$-periodic. As before, its fingerprint plot reveals that the dominating points belong to the set 
\begin{equation}
\label{eq:Am}
A_M = \{1\} \cup \{nM\pm 1 \, | \, n \in \mathbb{N} \},
\end{equation}
as a result, we compare the curve $z_M$, i.e., $z(t)$ rotated by $\pi/2-\pi/M$ in the counter-clockwise direction, with 
\begin{equation}
\label{eq:phiM}
\phi_M(t) = \sum_{k\in A_M} \frac{e^{2\pi i k^2 t}}{k^2}, \ t\in[0,1].
\end{equation}
\begin{figure}[htbp!] 
\centering
\includegraphics[width=0.456\textwidth,clip=true, valign=c]{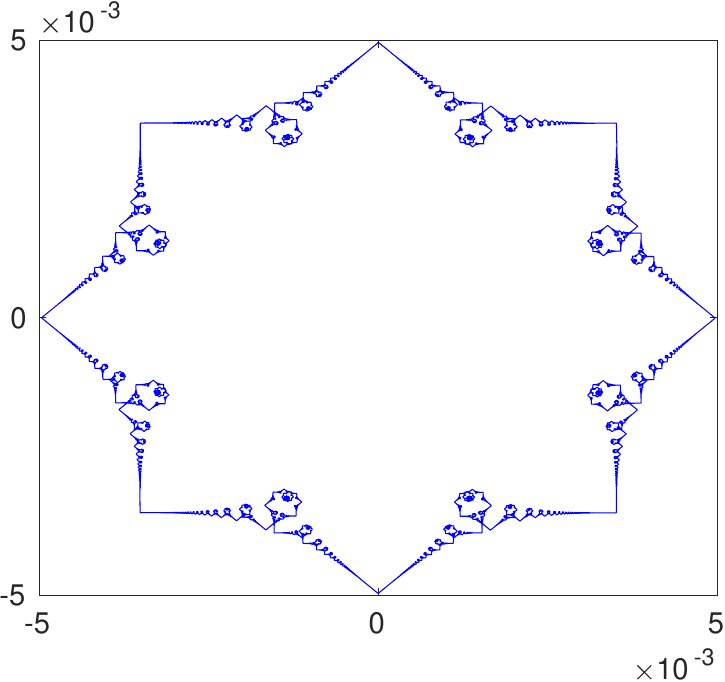}
\includegraphics[width=0.477\textwidth,clip=true, valign=c]{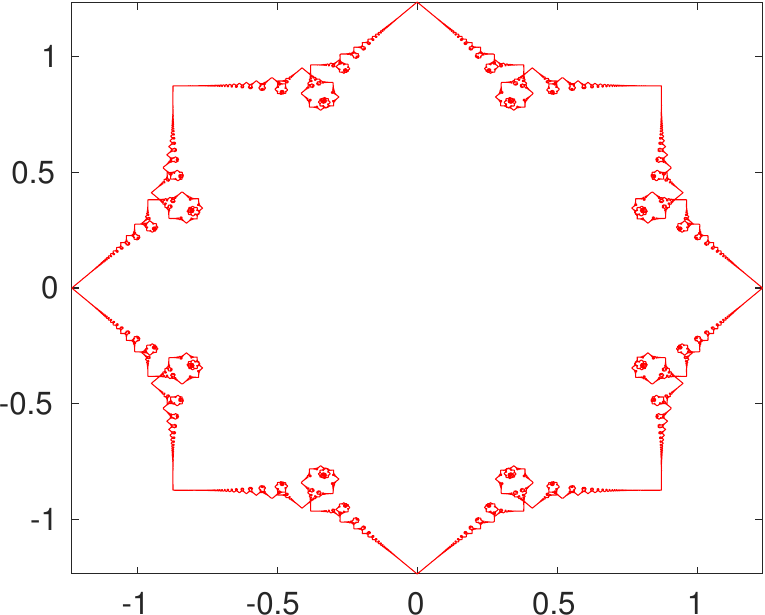}
\caption[]{Left: Trajectory of a single point $z_M(t)$, $M=4$, $N/M=2^{10}$. Right: $\phi_M(t)$ as in \eqref{eq:phiM} with $\dim(A_M)=2^{10}$. Both are computed at $N_t+1$ points in their respective domain intervals.}
\label{fig:zM3b1}   
\end{figure}
\begin{figure}[!htbp]\centering
\includegraphics[width=0.525\textwidth,clip=true, valign=b]{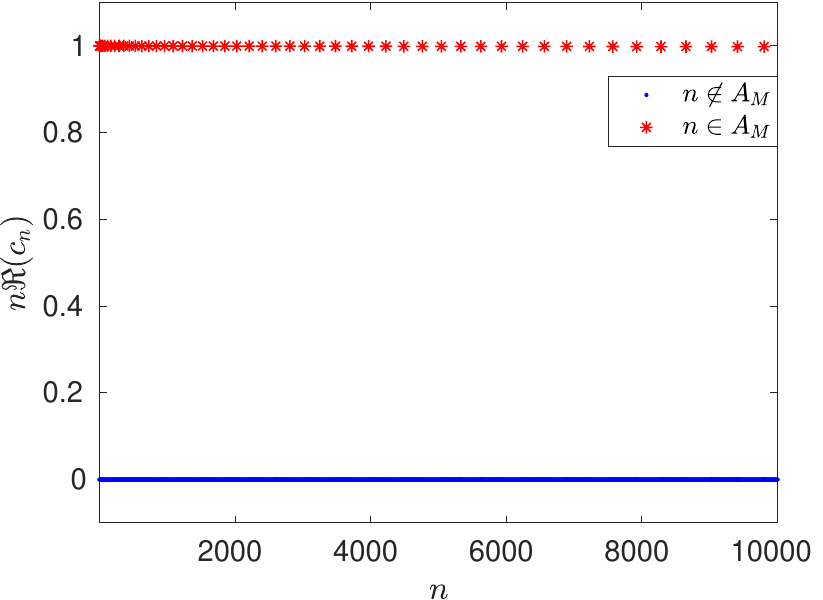}
\includegraphics[width=0.467\textwidth, clip=true, valign=b]{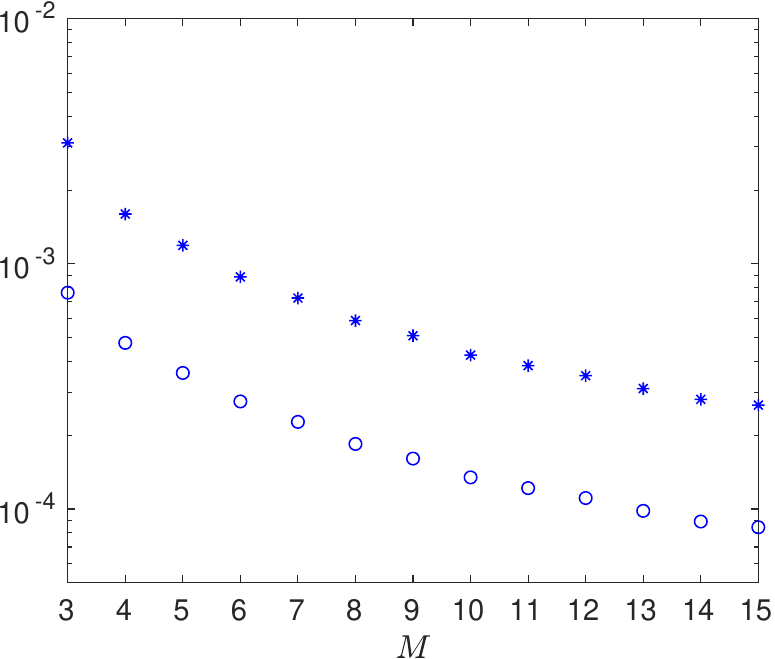}
\caption{Left: the fingerprint for $z_M$, $M=4$. Right: errors $\max_t|(\phi_M - \lambda_M  z_M - \mu_M)|$ (circled) and $\max_t|(\phi_M - \lambda_M  z_M - \mu_M)/ \phi_M|$ (starred), where $\lambda_M$ and $\mu_M$ are computed from \cite[(73)]{HozVega2014} using a least-squares fitting for $M=3,4,\ldots,15$.}
\label{fig:X0-phiM-errors-hyp}
\end{figure}
For $M=4$, on the left-hand side of Figure \ref{fig:zM3b1}, we have plotted $z_M(t)$ for $b=1+10^{-5}$, $N/M=2^{10}$, and $\phi_M$ on the right-hand side. One can clearly see that up to a scaling the two figures are very similar. The case $M=4$ is particularly interesting as the set $A_M$ involves only odd integers. To highlight further the comparison between $z_M$ and $\phi_M$, for different values of $M$, in Figure \ref{fig:X0-phiM-errors-hyp} (the right-hand side), we plot the absolute error $\max_t|(\phi_M - \lambda_M  z_M - \mu_M)|$ (circled) and relative error $\max_t|(\phi_M - \lambda_M  z_M - \mu_M)/ \phi_M|$ (starred) for $M=3, 4, \ldots, 15$. Here, $\lambda_M$ and $\mu_M$ are computed from \cite[(73) ]{HozVega2014} using least-square fitting. To make a fair comparison we have kept $N/M=2^{10}$ for all values of $M$ and since the step size $\Delta t = \mathcal{O}(1/M^2)$, the best agreement is expected for the larger values of $M$ as confirmed in our results which also shows the agreement between the two curves.  

The left-hand side of Figure \ref{fig:X0-phiM-errors-hyp} shows the fingerprint of the scaled $z_M$ where the dominating points in the set $A_M$ are very close to the value 1 and the rest close to the value 0. Very similar observations have been made for different values of $M$ as well which allows to conjecture that 
\begin{equation}
\lim\limits_{b\to 1^+} | n \, c_n|  = 
\begin{cases}
1 \ & \text{if $n \in A_M$,} \\
0 \ & \text{otherwise,} 
\end{cases}
\end{equation}
where $c_n$ are the Fourier coefficients of the scaled $z_M$; in other words, as $b\to1^+$, $z_M \to \phi_M$.

%
%
\section{$\T(s,\tpq)$ for $q\gg1$}
\label{sec:Tstpq-qgg1}
In Section \ref{sec:Pb-def-form}, evaluating \eqref{eq:psi-theta-st-hyp-hel} at times $\tpq$, i.e., the rational multiples of $T_f$ through analytical-algebraic techniques allowed us to obtain an expression involving generalized quadratic Gau{\ss} sums, and thus the algebraic solution. This would not be possible in case of any irrational time; however, it can be approximated with an arbitrary precision by some $\tpq$ with a large $q$. In this direction, we choose $\tpq$ with a \textit{large} $q$ such that there is no pair $\tilde{p}$, $\tilde{q}$, with both $\tilde{q}$ and $|\tfrac{p}{q} - \tfrac{\tilde{p}}{\tilde{q}}|$ \textit{small} and compute the evolution algebraically as it is free from numerical errors \cite{HozVega2014}.  

For the CHP, we take $M=3$, $\tpq = T_f(1/4 + 1/41 + 1/401)$, and after taking different values of the parameter $b$, we note that as $b$ tends to the value 1, the tangent vector $\T$ values move towards the south pole of $\Hbb^2$. The left-hand side of Figure \ref{fig:Tirr-hyp-hel} shows the stereographic projection of $\T$ about the $x$-axis for $b=1.2$ and one can see that the $Mq/2 = 98646$ values form a triskekion-like shape. Moreover, the spiral-like structures at a smaller scale exhibit a fractal-like phenomenon.  A similar behaviour is noticed in the case of a HHP as shown on the right-hand side of the same figure where $\tpq = (T_f) (1/3 + 1/31 + 1/301)$, $M=8$, $b=0.4$, $l=0.6$. Numerical experiments show that as the parameter $b$ grows the corresponding stereographic projection of $\T$ becomes steeper as one approaches from the end points to the center (compare it with the case $b=0$ in \cite[Fig. 11]{HozKumarVega2020}).

\begin{figure}[!htbp]\centering
\includegraphics[width=0.25\textwidth,  valign=c, clip=true]{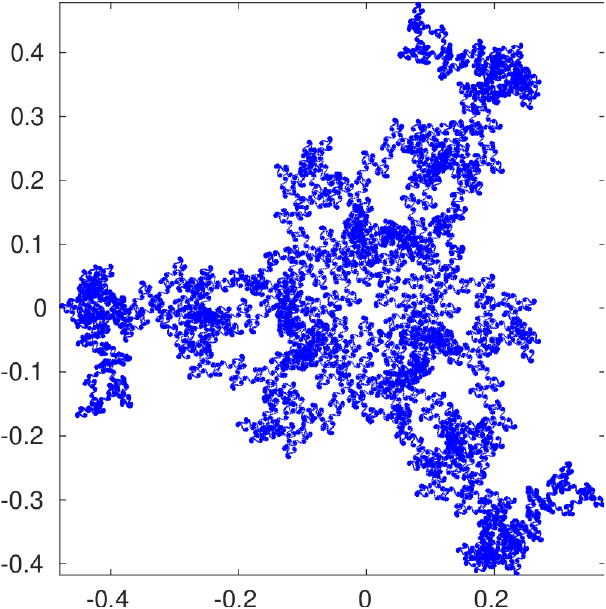}
\includegraphics[width=0.65\textwidth, valign=c, clip=true]{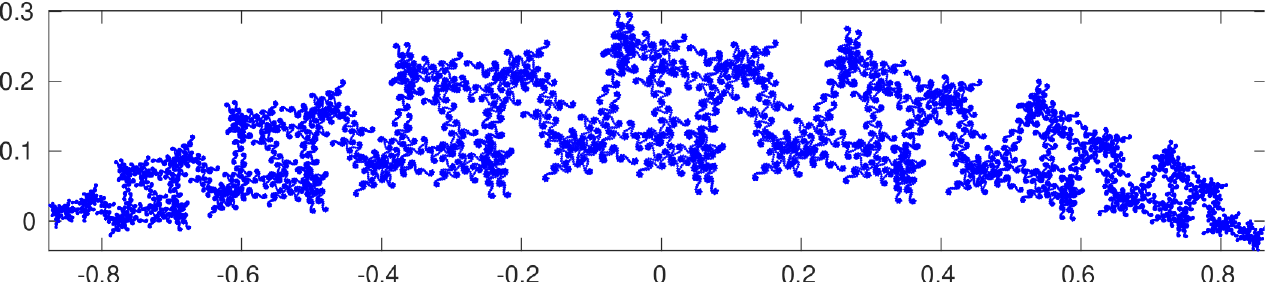}
\caption[]{The stereographic projection of tangent vector, i.e., $\left(\frac{T_2}{1+T_1},\frac{T_3}{1+T_1}\right)$. Left: $M=3$, $b=1.2$, $p=18209$, $q=65764$. Right: $M=8$, $b = 0.4$, $l=0.6$, $p=10327$, $q=27993$.}
\label{fig:Tirr-hyp-hel}
\end{figure}

\section{Numerical relationship between the helical polygon and one-corner problems}
\label{sec:Num-rel}
Numerical experiments in both Euclidean and hyperbolic settings have shown that the evolution of \eqref{eq:SMP-hyp}--\eqref{eq:VFE-hyp} with a regular polygon as an initial datum can be understood as a superposition of several one-corner problems \cite{HozVega2018,HozKumarVega2019,HozKumarVega2020}. This implies that for a polygonal curve, the asymptotic behaviour of the solution close to any of the corners is described through the self-similar solution corresponding to the angle of that corner \cite{delahoz2007}. Moreover, the relationship between the two problems has shown several profound implications, for instance, in \cite{HozKumarVega2020}, this fact was used to compute a precise expression for the speed of center of mass analytically and to recover the parameter $c_0$ in \eqref{eq:c0-ctheta0}. Although obtaining $c_l$, $c_M$ with analytic techniques is beyond the scope of this article, in what follows,  we establish this relationship numerically for both CHP and HHP and utilize it to obtain the parameter $c_{\theta,0}$ in \eqref{eq:c0-ctheta0}. We also use it to compute a quantity which in the planar polygon problem determined the angle of the planar curve $\X(0,t)$, and we compare it with its Euclidean counterpart as explained below. 
 
To compare the two problems, for a given $\theta_0$ (which also determines $c_{\theta,0})$ we solve the one-corner problem numerically for time $t=\tpq$, $q\gg1$ and denote the solution by $\T_\theta$, $\nn_\theta$, $\bb_\theta$, $\X_\theta$. Here, we choose $c_{\theta,0}$ as in \eqref{eq:c0-ctheta0} so that the asymptotes $\lim_{s\to-\infty} \T_\theta(s) = \A^-$, $\lim_{s\to +\infty} \T_\theta(s) = \A^+$ form an angle same as that made by the corner of the helical polygon $\X(s,0)$ at $s=0$. Next, we rotate the solution in such a way that it has the same orientation as that of the helical polygon problem. This is achieved through the rotation matrix $\M$ so that $\T_{rot} = \M \cdot \T_\theta$ with $\T_{rot}^\pm(s) = \lim_{s\to\pm\infty} \T_{rot}(s) $ satisfies $\T_{rot}^\pm(s) = (a \cosh(l/2),\pm a\sinh(l/2),b)^T$, and $\X_{rot} = \X_0 + \M \cdot \X_{\theta}$, where $\X_0 = (0, a (l/2)/ \sinh(l/2),0)^T$, for the HHP. Similarly, for the CHP, $\T_{rot}^-(s) = (b, a \cos(2\pi/M), -a \sin(2\pi/M))^T$, $\T_{rot}^+(s) = (b, a , 0)^T$ and $\X_0 = (0, -a \pi/M, -a\pi / (M\tan(\pi/M)))^T$. Here, $\X_0$ corresponds to the location of the corner of the helical polygon $\X(s,0)$ at $s=0$. 

We solve the one-corner problem numerically using a fourth-order Runge--Kutta method with $\Delta s = l/2q$ for the HHP and $\Delta s = \pi / Mq$ for the CHP at $t=t_{1,q}$ with $q=4000$. With the same spatial step size we compute the evolution of \eqref{eq:SMP-hyp}--\eqref{eq:VFE-hyp} for the CHP with $M=6$, $b=1.2$, and with $M=8$, $l=0.6$, $b=0.4$ for the HHP. We also compare the curve $\X$, in particular, the evolution $\X(0,t)$ for $t\in[0,t_{1,20}]$. Figure \ref{fig:M-1-corner-circ-hel}  shows the CHP case where the left-hand side shows the stereographic projection of $\T$, and the right-hand side shows the projection of $\X(0,t)$ onto the complex plane. The blue color corresponds to the helical polygon problem and the red to the one-corner problem. The HHP case has been displayed in Figure \ref{fig:M-1-corner-hyp-hel}. One can clearly notice that for the tangent vector $\T$, except the center for spiral curves, the two curves match very well. Similarly, for small times, the projection of $\X(0,t)$ (in blue) can be well approximated with a straight line (in red) with a slope $A_2/A_3$. Recall that in \cite{HozKumarVega2020}, the quantity $A_2/A_3$, which depends only on $c_{\theta,0}$ through \cite[(74)--(75)]{HozKumarVega2020} determined the angle of $\X(0,t)$ which it makes with the plane containing $\X(s,0)$. As $A_2/A_3$ quantifies how $\X(0,t)$ grows near $t=0$, we use it to compare this growth for the helical polygon problem in both Euclidean and hyperbolic cases. We plot it for different values of $c_{\theta,0}$, for the two cases and the bottom right of Figure \ref{fig:M-1-corner-hyp-hel} shows that $A_2 /A_3$ is greater (respectively, smaller) than the one in the hyperbolic (respectively, Euclidean) case, and tends to one, as $c_{\theta,0}$ tends to zero.

Another important consequence of the relationship between the two problems is to recover the coefficient $c_0$ (now, $c_{\theta,0}$) in the one-corner problem from the multiple corner problem \cite{HozKumarVega2020}. Thus, for infinitesimal small times, from the one-corner problem, the curvature of the helical polygon $\X$ at $s=0$, can be written as $c_{\theta,0}(\tpq) = \sqrt{\tpq} |\Ts(0,\tpq)|_0$, where $\tpq\ll1$ \cite{delahoz2007}. We approximate the derivative using a second-order finite difference and note that, given any $q$, as $\theta_0>0$, the Galilean shift $s_{1,q} < 2\pi/Mq$ for the CHP and $s_{1,q}<l/q$ for the HHP. As a result, the tangent vector $\T(s,t_{1,q})$ is continuous at $s=0$, $s=-\Delta s$ and $s=\Delta s$, where $\Delta s = 2\pi/Mq \ (4\pi/Mq)$ if $q$ is odd (even) for the CHP and $\Delta s = l/q \ (l/2q)$ if $q$ is odd (even) for the HHP. This implies that $\T(s_{1,q}^-, t_{1,q}) = \T(0, t_{1,q})$, $\T((-\Delta s + s_{1,q})^-, t_{1,q}) = \T(-\Delta s, t_{1,q})$ and $\T((\Delta s + s_{1,q})^-, t_{1,q}) = \T(\Delta s, t_{1,q})$. Using the algebraic solution $\Talg$, for $q\equiv0\bmod4$, we can write 
\begin{equation}
\label{eq:c_theta_q-approx}
c_{\theta,0} \approx c_{\theta,0}^q \equiv \frac{|\Talg\left( \Delta s,t_{1,q}\right) -\Talg\left( - \Delta s,t_{1,q}\right)|_0 }{2 \Delta s}, \ q\gg1.
\end{equation}
We have computed $c_{\theta,0}^q$ numerically for the same parameter values as before with $q = 1000\cdot 2^{r}$, $r=0,1,2,\ldots,7$, and the absolute errors using \eqref{eq:c0-ctheta0} for both the cases. Table \ref{table:c_0t_error} clearly shows that the errors in both cases decrease by half as $q$ is doubled, thereby showing the behaviour $\mathcal{O}(1/q)=\mathcal{O}(t_{1,q})$, i.e., a first-order convergence. 

 \begin{table}[!htbp]
 \centering
 \begin{tabular}{|c|c|c|}
 \cline{2-3}
 \multicolumn{1}{c|}{} & \multicolumn{2}{c|}{$|c_{\theta,0} - c_{\theta,0}^q|$}
 \\
 \hline
 $q$ & CHP & HHP
 \\
 \hline
 1000 & $2.8194\cdot10^{-5}$ & $2.6883\cdot10^{-5}$
 \\
 2000 & $1.4181\cdot10^{-5}$ & $1.3630\cdot10^{-5}$
 \\
 4000 & $7.1117\cdot10^{-5}$ & $6.8621\cdot10^{-5}$
 \\
 8000 & $3.5611\cdot10^{-6}$ & $3.4428\cdot10^{-6}$
 \\
 \hline  
 \end{tabular}~\begin{tabular}{|c|c|c|}
 \cline{2-3}
 \multicolumn{1}{c|}{} & \multicolumn{2}{c|}{$|c_{\theta,0} - c_{\theta,0}^q|$}
 \\
 \hline
 $q$ & CHP & HHP
 \\
 \hline
 16000 & $1.7818\cdot10^{-6}$ & $1.7244\cdot10^{-6}$
 \\
 32000 & $8.9134\cdot10^{-7}$ & $8.6292\cdot10^{-7}$
 \\
 64000 & $4.4590\cdot10^{-7}$ & $4.3164\cdot10^{-7}$
 \\
 128000 & $2.2204\cdot10^{-7}$ & $2.1586\cdot10^{-7}$
 \\
 \hline  
 \end{tabular}	
 
 \caption{$|c_{\theta,0} - c_{\theta,0}^q|$, for a CHP for $b=1.2$, $M=6$, $c_{\theta,0}=0.1823\ldots$ and for a HHP for $b=0.4$, $M=8$, $l=0.6$, $c_{\theta,0}=0.1803\ldots$ for different values of $q$, where $c_{\theta,0}^q$ is given by \eqref{eq:c_theta_q-approx}. In both cases, the errors decrease as $\mathcal{O}(1/q) = \mathcal{O}(t_{1,q}).$}
 \label{table:c_0t_error}
 \end{table}
	
\begin{figure}[!htbp]\centering
\includegraphics[width=0.405\textwidth, valign=t, clip=true]{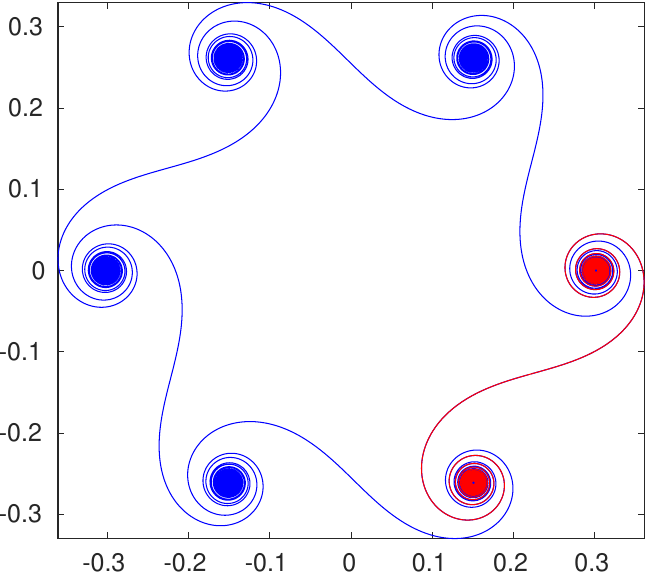}
\includegraphics[width=0.45\textwidth,  valign=t, clip=true]{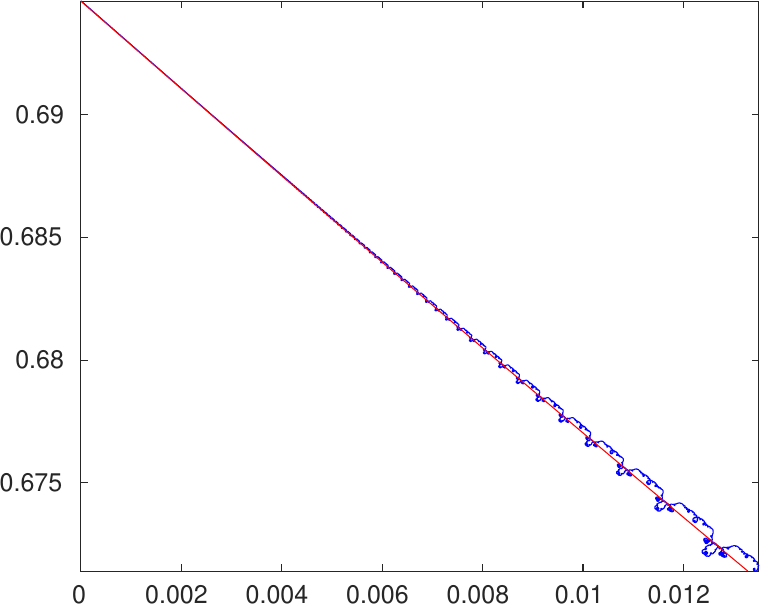}
\caption[]{For a CHP with $M=6$, $b=1.2$, $c_0=0.2082\ldots$; left: The stereographic projection of tangent vector at $t=t_{1,4000}$. Right: $-X_1(0,t)+i \|X_2(0,t)+X_3(0,t)\|_2$, for $t\in[0,t_{1,20}]$. }
\label{fig:M-1-corner-circ-hel}
\end{figure}
\begin{figure}[!htbp]\centering
\includegraphics[width=\textwidth, valign=t, clip=true]{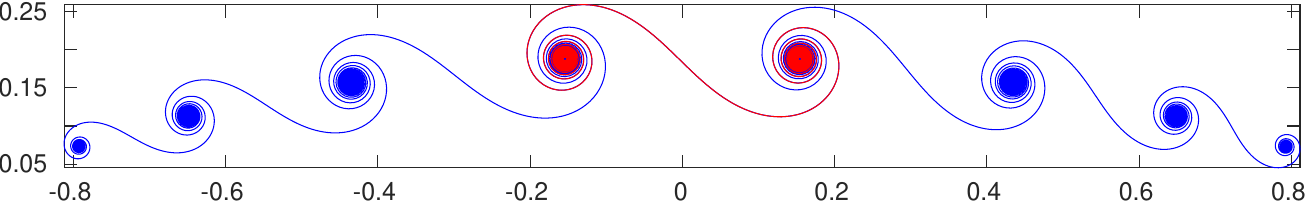}\\\vspace{0.2cm}
\includegraphics[width=0.45\textwidth,  valign=t, clip=true]{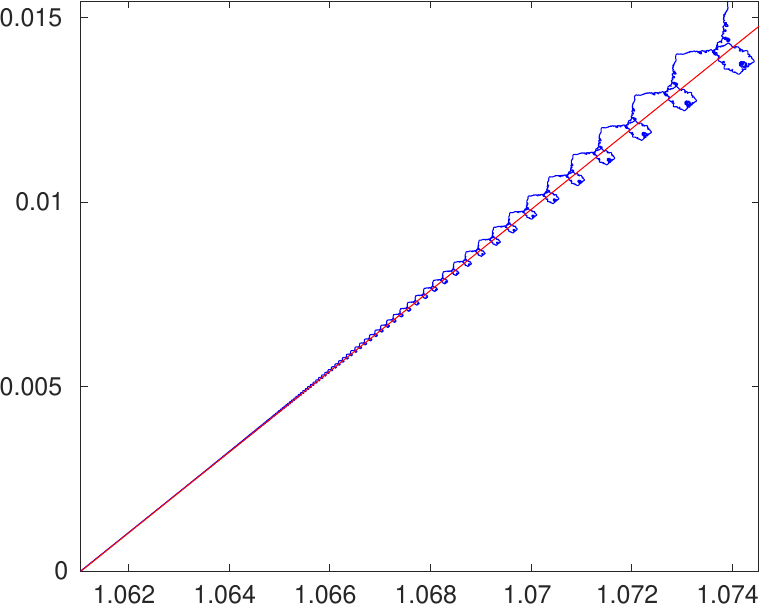}
\includegraphics[width=0.52\textwidth,  valign=t, clip=true]{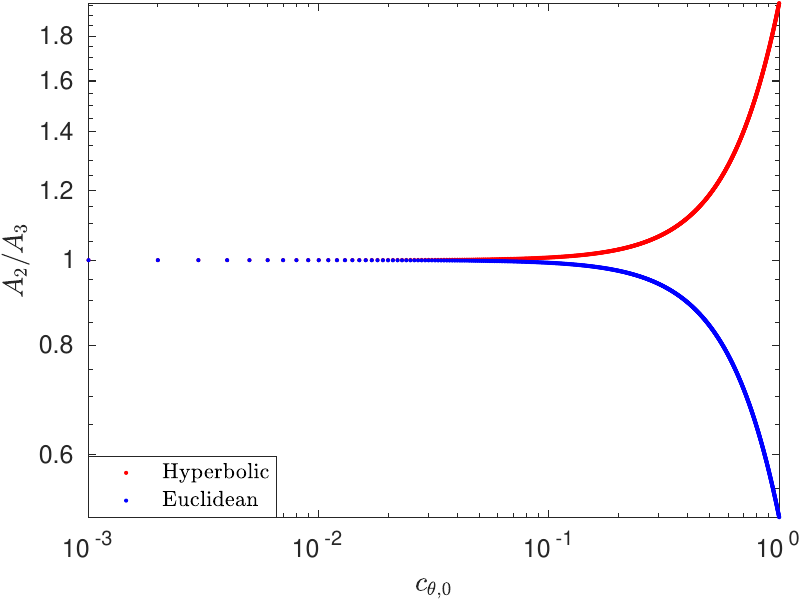}
\caption[]{For a HHP with $M=8$, $b = 0.4$, $c_0=0.1803\ldots$; top: the stereographic projection of tangent vector at $t=t_{1,4000}$. Left: $|X_1(0,t),X_2(0,t)|_0+iX_3(0,t)$. Right: A log-log plot of $A_2/A_3$ as a function of $c_{\theta,0}$ for both Euclidean and hyperbolic cases, which converges to one as $c_{\theta,0}$ goes to zero.}
\label{fig:M-1-corner-hyp-hel}
\end{figure}
\section*{Acknowledgments}
The author would like to thank Luis Vega, Carlos J. Garc\'{i}a-Cervera and Ananyo Dan for their advice and support, also to the anonymous referees whose suggestions have greatly improved this article. This work has been supported by the ERCEA Advanced Grant 2014 669689 - HADE, Severo Ochoa grant SEV-2017-0718, and the Basque Government BERC Program 2018-2021.
\bibliography{ArticleRevisedCleanFinal}
\bibliographystyle{ieeetr}

\end{document}